%% file: LanczosCoeffs.tex
\documentclass[a4paper,onecolumn,12pt,epsf]{article}
\usepackage{a4}
\usepackage{array}
\usepackage{epsfig}
\usepackage{epstopdf}
\usepackage{rotating}
\usepackage{amsmath}
\usepackage{amsfonts}
\usepackage{longtable}
\usepackage{url}
\usepackage{natbib}         
\begin{document}

\title{The Lanczos Approximation for the $\Gamma$-Function with Complex Coefficients}
\author{ William Rea \\
Christchurch, New Zealand \\
New Zealand  \\
email: rea.william@gmail.com}

\maketitle

\input Abstract.tex

\input Introduction.tex

\input LanczosCoeffsProps.tex

\input Numerical.tex

\input Conclusions.tex

\bibliography{gamma}

\clearpage

\appendix

\input AllGraphs.tex


\end{document}

%% file: Abstract.tex
\begin{abstract}
We examined the properties of the coefficients of the \cite{lanczos1964}
approximation of the $\Gamma$-function with complex value of the free
parameter together with the convergence properties of the approximation
when using these coefficients. We report that for fixed real parts of the
free parameter that using complex coefficients both increases the 
computational cost of the Lanczos approximation while drecreasing the
accuracy. We conclude that in practical applications of numerical
evaluation of the $\Gamma$-function only coefficients generated with
real values of the free parameter should be used.. 
\end{abstract}

%% file: Introduction.tex
\section{Introduction}

The $\Gamma$-function, defined by
\begin{equation}
\Gamma(z)=\int_0^{\infty}t^{z-1}e^{-t}dt\label{eqn:maindef}
\end{equation}
was introduced by Euler as a solution to the ``interpolation''
problem of defining the factorial function for non-integer values.
Equation \eqref{eqn:maindef} has the desired property that $\Gamma(n+1)=n!$. 
The uniqueness of this definition was proved by \cite{bohr1922} and
a second proof attributed to H. Wielandt is discussed by \cite{Remmert1996}
and included in the textbook by \cite{Freitag2005},

In
the form of Equation \eqref{eqn:maindef} it is defined in the complex plane 
for $\mathcal{R}e(z)>0$
and can be extended to the left half-plane, among other several 
other methods, through the use of the
Euler reflection formula
\begin{equation}
\Gamma(z)\Gamma(1-z)=\frac{\pi}{\sin\pi z}.\label{eqn:reflection}
\index{sine}
\end{equation}

The $\Gamma$-function has found wide application in pure and applied
mathematics, statistics and the physical sciences.
Exact values for the $\Gamma$-function can be found for positive 
integer (zero and the negative integers are simple poles) and for
half integer values of $z$. Apart from those values, the $\Gamma$-function
must be evaluated using numerical methods. To this end a number of 
approximations and asymptotic series have been developed. The earliest
being Stirling's approximation; for integer $n$
$$
n!\approx \sqrt{2\pi n}n^ne^{-n}=\sqrt{2\pi n}\left(\frac{n}{e}\right)^n.
$$
A large number of proofs of Stirling's approximation and
closely related formulas are available, see, for example,
\cite{Aissen1954},
\cite{Feller1967},
\cite{Spira1971},
\cite{Namias1986},
\cite{Marsaglia1990},
\cite{Deeba1991},
\cite{Schuster2001},
\cite{Michel2002},
\cite{Dutkay2013},
\cite{Lou2014}, and
\cite{Neuschel2014}.

There are also a number of series expansions available such as 
Stirling's asymptotic series 
$$
n!=\left(\frac{n}{e}\right)^n\sqrt{2\pi n}\exp\left(\frac{1}{12n}
-\frac{1}{360n^3}+\frac{1}{1260n^5}-\dots\right).
$$
There are some proofs
which only require elementary methods such as that of \cite{Mermin1984},
\cite{Namias1986}, \cite{diaconis1986}, \cite{Patin1989}
and \cite{Marsaglia1990}.
 
In addition there are series methods due to  \cite{Spira1971},
 \cite{Spouge1994}, and \cite{Schmelzer2007}. Among these other series
\cite{lanczos1964} derived an expansion for the
numerical evaluation of the
$\Gamma$-function which was popularized to some extent by its inclusion
in \cite{Press1997} and which was studied in some detail by \cite{Pugh2004}.
The expansion has two notable features; it is convergent rather than an
asymptotic expansion and it has a free parameter (denoted by
$\gamma$ by \cite{lanczos1964} but by $r$ by \cite{Pugh2004}, we will use
$r$ in this paper) which affects the
rate of convergence of the approximation. While it is known that the
free parameter can be complex, to the best of this author's knowledge
the convergence properties of a complex expansion have not been
studied. \cite{Godfrey2001} studied the convergence properties of
the expansion for real values of $r$ and implemented his best 
choice of coefficients in Matlab \citep{MATLAB} using
15 terms and giving an approximation that was
accurate to 15 significant digits along the real axis
and 13 significant digits elsewhere.

The Lanczos expansion,, which is for $\Gamma(z+1)$ rather than $\Gamma(z)$, 
can be written
\begin{multline}
\Gamma(z+1)=
\sqrt{2\pi}\left(z+r+\frac{1}{2}\right)^{z+\frac{1}{2}}
             e^{-\left(z+r+\frac{1}{2}\right)} \\
    \left[\frac{a_0(r)}{2}+\frac{z}{z+1}a_1(r)+\frac{(z-1)z}{(z+2)(z+1)}a_2(r)
   +\cdots\right]
\label{eqn:LanczosSeries}
\end{multline}
where the $a_k(r)$ are the coefficients 
and $r$ is the free parameter. 

The key
to obtaining the $a_k(r)$ is to note that the $\Gamma$-function has known
values at the integers and that for positive integers the expansion 
\eqref{eqn:LanczosSeries}
terminates, allowing the $a_k(r)$ to be determined recursively.
\cite{Pugh2004} investigated several methods of obtaining
the $a_k(r)$ and reported the recursive methods worked well in practice.

In \eqref{eqn:LanczosSeries} if $z=0$ then $\Gamma(1)=1$ and the
series terminates after the first term. Thus
\begin{align}
1&=\sqrt{2\pi}\left(r+\frac{1}{2}\right)^{1/2}
  e^{-\left(r+\frac{1}{2}\right)}\frac{a_0(r)}{2}.\nonumber\\
\intertext{Solving for $a_0(r)$ we obtain}
a_0(r)&=\sqrt{\frac{2e}{\pi(r+\frac{1}{2})}}e^r.\label{eqn:CalcA0r}
\end{align}
Having obtained $a_0(r)$ we then set $z=1$ so that $\Gamma(2)=1$ and
\eqref{eqn:LanczosSeries} now terminates after the second term allowing
us to find $a_1(r)$. By successively setting $z=2,3,4\ldots$ we can
calculate as many of the $a_k(r)$ as required.

On a practical
matter of finding the $a_k(r)$ it is useful to set
\begin{align}
F_r(z)&=\frac{\Gamma(z+1)e^{z+r+1/2}}{\sqrt{2\pi}
             \left(z+r+1/2\right)^{z+1/2}}.\label{eqn:LanczosFrz}
\intertext{Then}
a_1&=\left(F_r(1)-\frac{a_0}{2}\right)\frac{2}{1} \label{eqn:LanczosA1} \\
a_2&=\left(\left(F_r(2)-\frac{a_0}{2}\right)\frac{3}{2}-a_1\right)\frac{4}{1}
\label{eqn:LanczosA2}
\intertext{The general ($n\ge 3$) recursion formula is}
a_n(r)&=\left(\left(\left(\left(F_r(n)-\frac{a_0}{2}\right)\frac{n+1}{n}
       -a_1\right)\frac{n+2}{n-1}
       -a_2\right)\frac{n+3}{n-2} - \cdots - a_{n-1}\right)\frac{2n}{1}.
\label{eqn:LanczosAn}
\end{align}
The final equation is \cite{Pugh2004} Equation (6.4).

Having found some $a_k(r)$ values they are substituted
into the  expansion \eqref{eqn:LanczosSeries} and desired approximations
of $\Gamma(z+1)$ can be obtained, values with $\mathcal{R}e(z)<0$ can
be obtained by combining the Lanczos approximation with the
Euler reflection formula, Equation \eqref{eqn:reflection}.

In the remainder of this paper Section (\ref{sec:Properties}) investigates
some of the properties of the Lanczos coefficients with complex
values of the free parameter, Section (\ref{sec:Numerical}) 
reports on the numerical
performance of series with complex coefficients, Section (\ref{sec:Conclude})
presents our conclusions while Appendix (\ref{sec:Supplement}) presents
some supplementary materials.

%% file: LanczosCoeffsProps.tex
\section{Properties of the Complex $a_k(r)$}\label{sec:Properties}

The formulas \eqref{eqn:CalcA0r} through \eqref{eqn:LanczosAn} allow
us to calculate the expansion coefficients with complex $r$.
All numerical values and graphics
presented here were obtained with custom written
Matlab \citep{MATLAB} code. The code was verified using published values
of the $a_k(r)$ for real $r$ in \cite{lanczos1964} and
\cite{Pugh2004} and then, without any further alteration of
the code, $r$ was allowed to be complex and the required 
coefficients calculated.

\begin{figure}[h]
\begin{center}
\includegraphics[width=12cm]{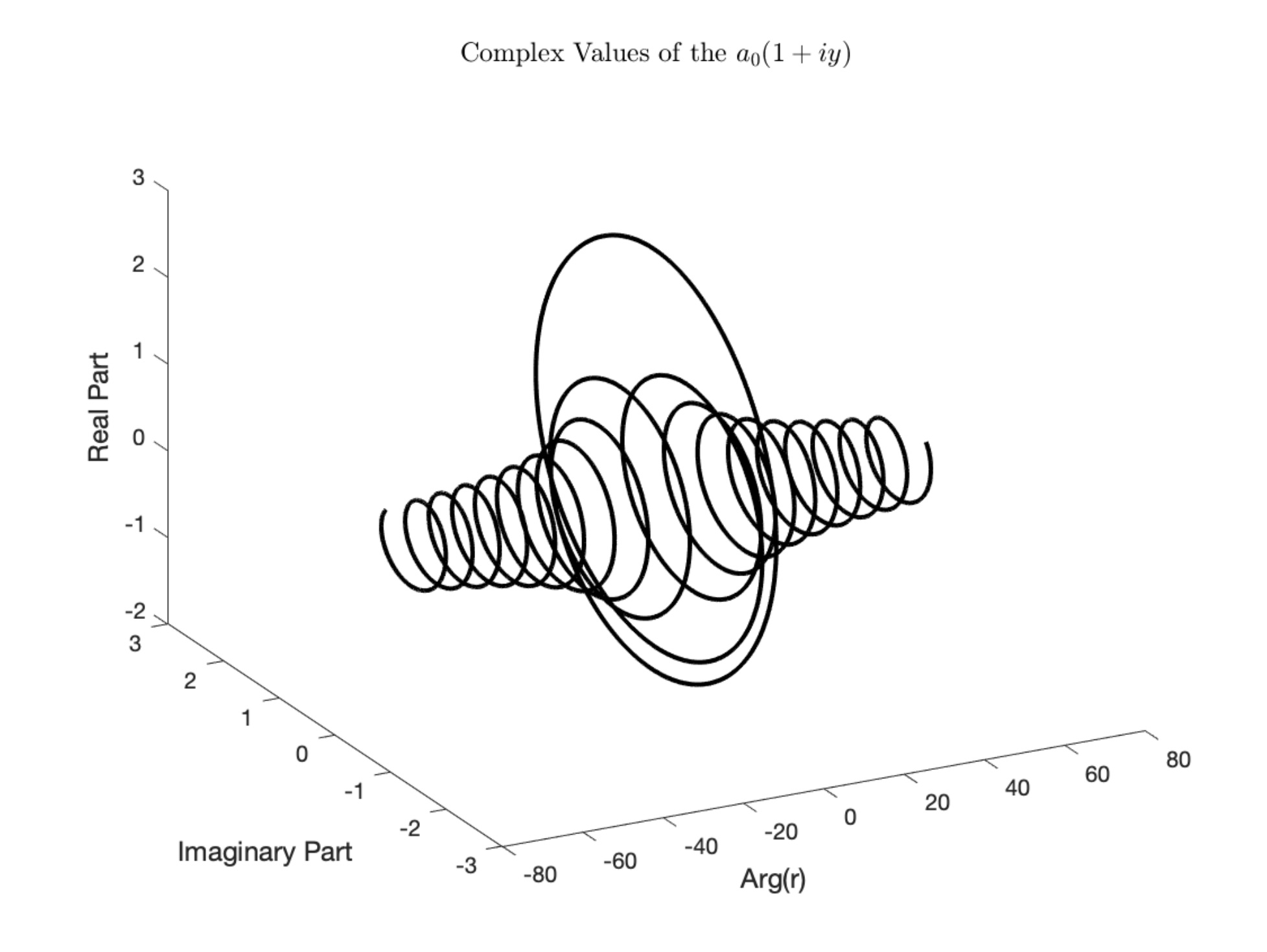}
\end{center}
\caption{The complex values of the $a_0(r)$ coefficient of the
Lanczos expantion with $r=1+ir_y$ where $-20\pi\le r_y \le 20\pi$.
}\label{fig:LanczosA0Complex}
\end{figure}

\begin{figure}[h]
\begin{center}
\includegraphics[width=12cm]{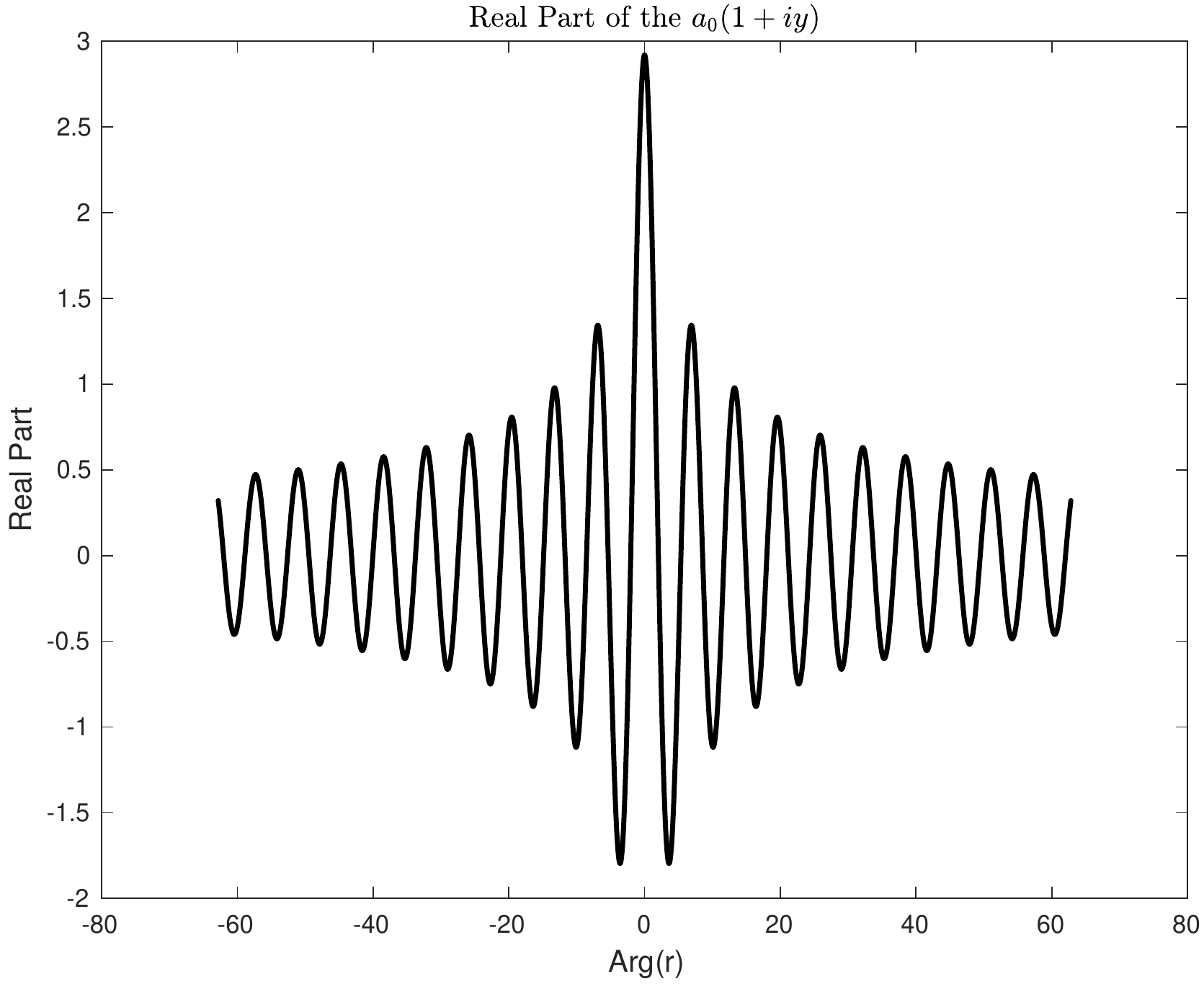}
\end{center}
\caption{The real values of the $a_0(r)$ coefficient of the
Lanczos expantion with $r=1+ir_y$ where $-20\pi\le r_y \le 20\pi$.
}\label{fig:LanczosA0Real}
\end{figure}

\begin{figure}[h]
\begin{center}
\includegraphics[width=12cm]{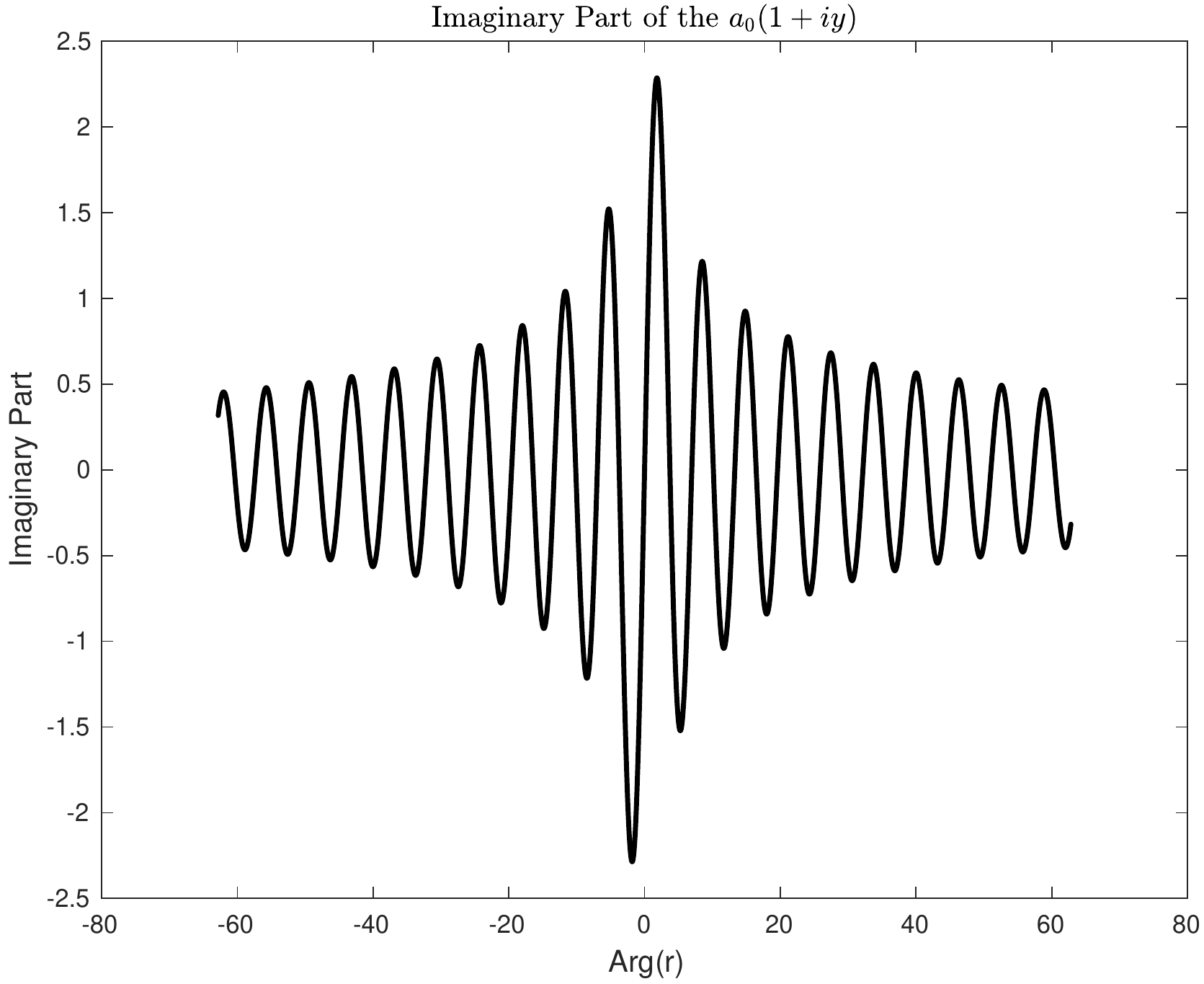}
\end{center}
\caption{The imaginary values of the $a_0(r)$ coefficient of the
Lanczos expantion with $r=1+ir_y$ where  $-20\pi \le r_y \le 20\pi$.
}\label{fig:LanczosA0Imag}
\end{figure}

\begin{figure}[h]
\begin{center}
\includegraphics[width=12cm]{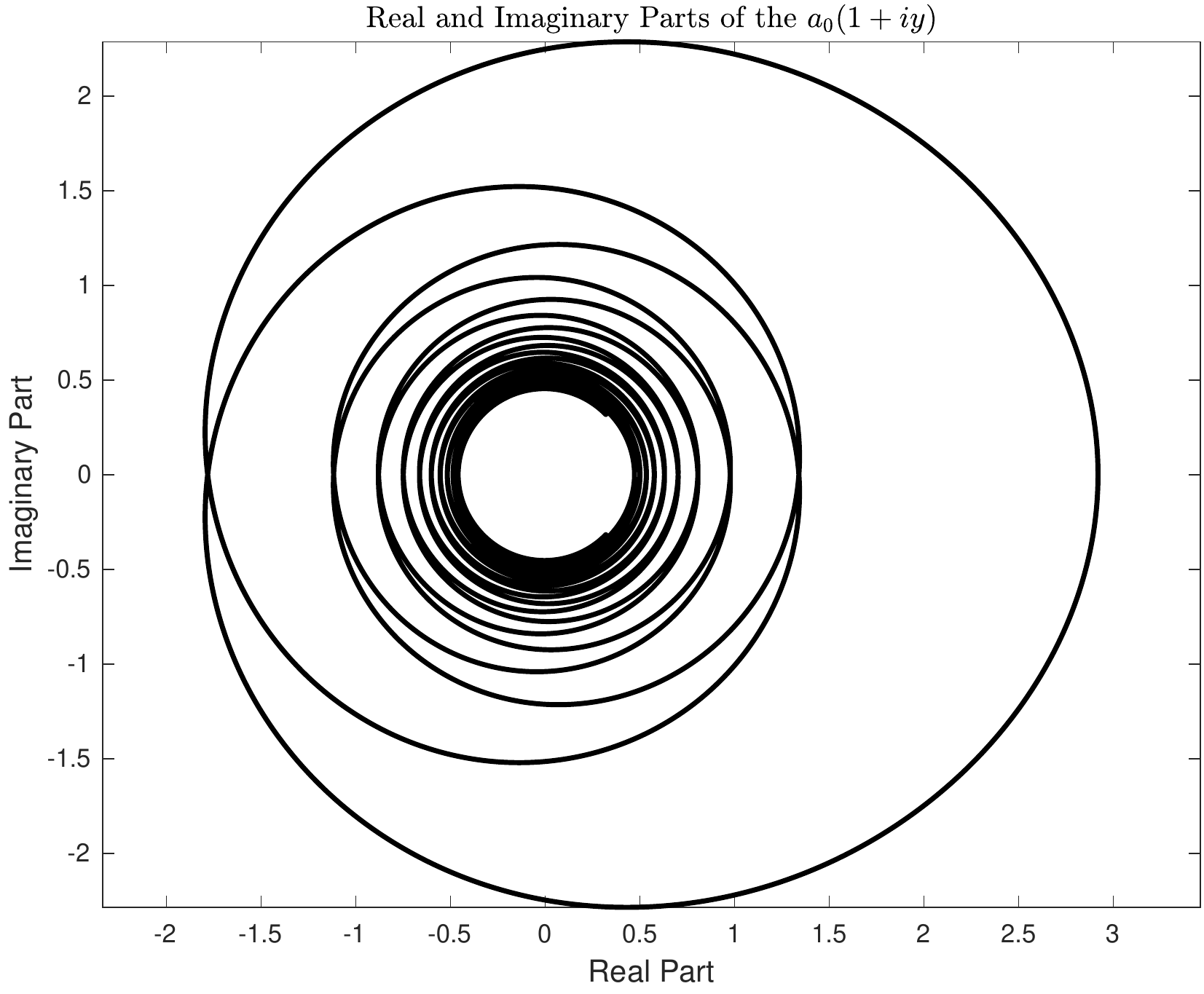}
\end{center}
\caption{The real and imaginary values of the $a_0(r)$ coefficient of the
Lanczos expantion with $r=1+ir_y$ where  $-20\pi \le r_y \le 20\pi$.
}\label{fig:LanczosA0ReIm}
\end{figure}

\begin{figure}[h]
\begin{center}
\includegraphics[width=12cm]{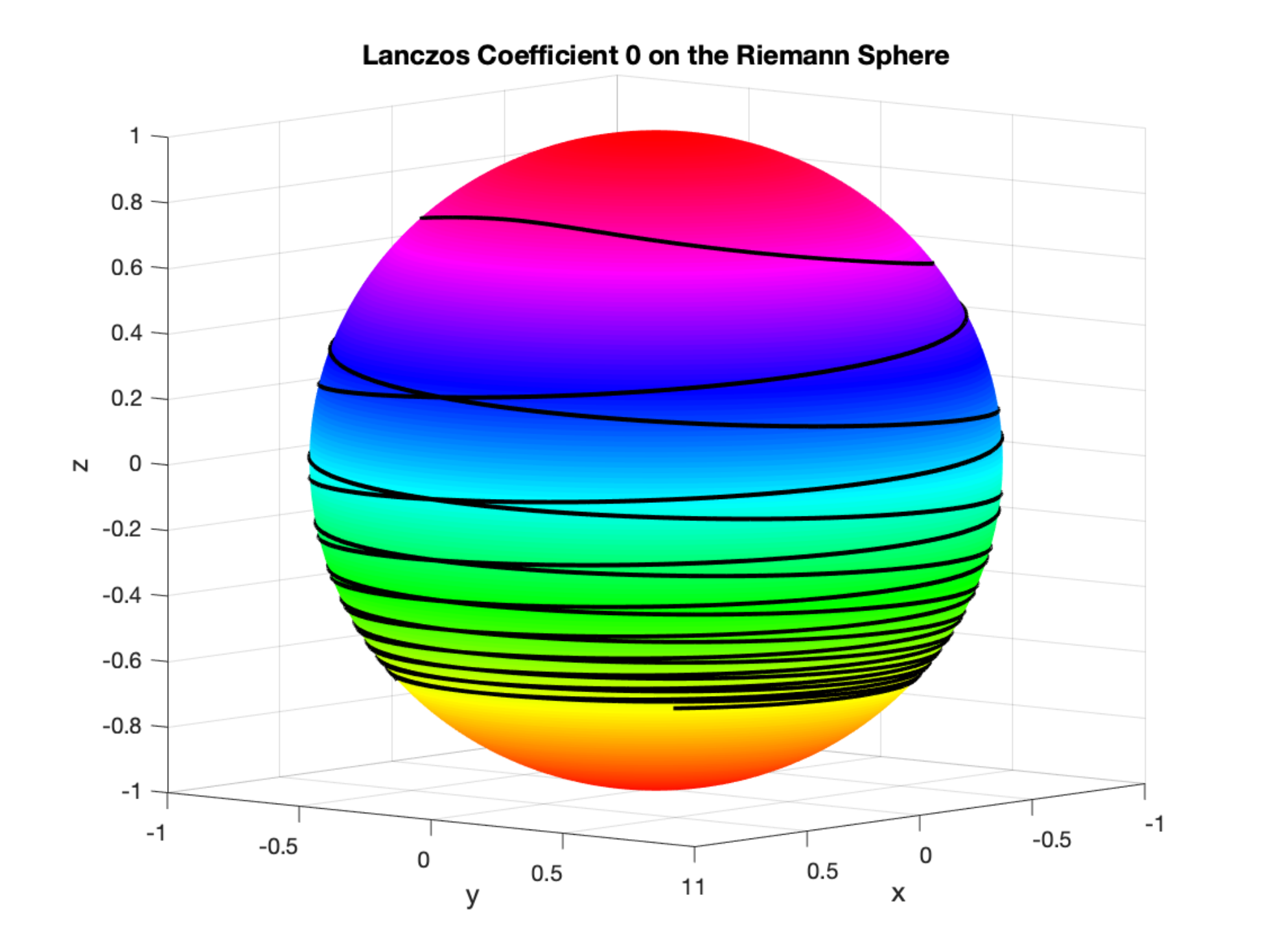}
\end{center}
\caption{The coefficient
$a_0(r)$ for $r=1+ir_y$ for $-20\pi \le r_y\le 20\pi$ plotted on the
Riemann sphere.}\label{fig:A0RS}
\end{figure}

Figures (\ref{fig:LanczosA0Complex}) through (\ref{fig:A0RS}) present
five different views of the $a_0(r)$ coefficient for $r=1+ir_y$ with
$-20\pi \le r_y \le 20\pi$. While we explored values of $r$ with different
real parts ($\mathcal{R}e(r)=1,3,\ldots,13$), the behaviour of the 
coefficients were similar in all cases,
the main change being their magnitudes.

In the recursive calculation of the $a_k(r)$ only integer values of $z$ are
used so we denote these by $n=1,2,3,\ldots$.  The free parameter
enters at two points. In the numerator of the $F_r(z)$, 
Equation \eqref{eqn:LanczosFrz}, with 
$n=1,2,3,\ldots$, has the term $e^{n+r+1/2}$. Writing $r=r_x+ir_y$ the
real part of the exponent is $n+r_x+1/2$ while the imaginary
part is $ir_y$. Holding the real part constant and allowing $\mathcal{I}m(r)$
to change makes the coefficients ``spin'' in the complex plane, examples
can be seen in Figures (\ref{fig:LanczosA0Complex}), 
(\ref{fig:LanczosA0ReIm}) and (\ref{fig:A0RS}).

The other place $r$ enters the calculation of the $a_k(r)$ is in the
denominator of $F_r(z)$. This is what causes the ``tappering off'' seen
in Figure (\ref{fig:LanczosA0Complex}) through (\ref{fig:LanczosA0Imag}) and the
gradual descent of the path of the $a_k(r)$ on the Riemann sphere, see Figure
(\ref{fig:A0RS}).
If
\begin{align*}
r&=r_x+ir_y
\intertext{We can rewrite this in polar form}
r&=Re^{i\theta}
\intertext{ where}
R&=\sqrt{r_x^2+r_y^2} \\
\theta&=\tan^{-1}\left(\frac{r_y}{r_x}\right).
\end{align*}
If $r_x$ is fixed increasing $|r_y|$ increases $R$ causing the tapper.
Asympotitically the real and imaginary parts of $a_0(r)$
decay as $1/\sqrt{r_y}$.

Figures (\ref{fig:LanczosA0Real}) and (\ref{fig:LanczosA0Imag}) suggest
that there are points for fixed $r_x$ ($r_x=1$ in the Figures) for
 which the $a_0(r)$ will be purely real or
purely imaginary. We can derive conditions under which these cases occur
and show that there are an infinite number of such coefficients.

Using Equation \eqref{eqn:CalcA0r}, if $a_0(r)$ is purely imaginary then
it has the form $a_0(r)=iv$. Thus
\begin{align*}
iv&=\sqrt{\frac{2e}{\pi(r+\frac{1}{2})}}e^r \\
-v^2&=\frac{2e}{\pi(r+\frac{1}{2})}e^{2r}.
\intertext{First, splitting $r$ into its real and imaginary parts we have}
-v^2&=\frac{2e}{\pi\left((r_x+\frac{1}{2})+ir_y\right)}e^{2r_x+i2r_y}.
\intertext{Next, we proceed to split this into its real and imaginary
parts}
-v^2&=\frac{2e^{2r_x+1}}{\pi\left(\left(r_x+\frac{1}{2}\right)^2+r_y^2\right)}
  \left(\left(r_x+\frac{1}{2}\right)-ir_y\right)
  \left(\cos 2r_y+i\sin 2r_y\right)
\end{align*}
Thus we have
\begin{multline*}-v^2= \frac{2e^{2r_x+1}}
                   {\pi\left(\left(r_x+\frac{1}{2}\right)^2+r_y^2\right)}
  \Bigg[\left(\left(r_x+\frac{1}{2}\right)\cos 2r_y+r_y\sin 2r_y\right) \\
   +i\left(\left(r_x+\frac{1}{2}\right)\sin 2r_y-r_y\cos 2r_y\right)\Bigg].
  \end{multline*}
From this we neglect the sign of the leading fraction because it is always
positive. This leaves us with two conditions for the $a_0(r)$ to be purely 
imaginary.
The first is that real part must be negative, that is 
\begin{align}
\left(r_x+\frac{1}{2}\right)\cos 2r_y+r_y\sin 2r_y&<0 \nonumber \\
\intertext{or}
\tan 2r_y <-\frac{r_x+1/2}{r_y}.\label{eqn:Ak0Imagl}
\intertext{Secondly, we require that the imaginary part is zero}
\left(r_x+\frac{1}{2}\right)\sin 2r_y-r_y\cos 2r_y&=0. \nonumber
\intertext{From this second condition we require}
\tan 2r_y&=\frac{r_y}{r_x+1/2}. \label{eqn:Ak0Imag2}
\intertext{For the $a_0(r)$ to be purely real the condition on the
imaginary part (Equation \ref{eqn:Ak0Imag2}) is unchanged and in place of
Equation \eqref{eqn:Ak0Imagl} we require}
\left(r_x+\frac{1}{2}\right)\cos 2r_y+r_y\sin 2r_y&>0. \nonumber
\intertext{or}
\tan 2r_y <-\frac{r_x+1/2}{r_y}.  \label{eqn:Ak0Real}
\end{align}

\begin{figure}[h]
\begin{center}
\includegraphics[width=12cm]{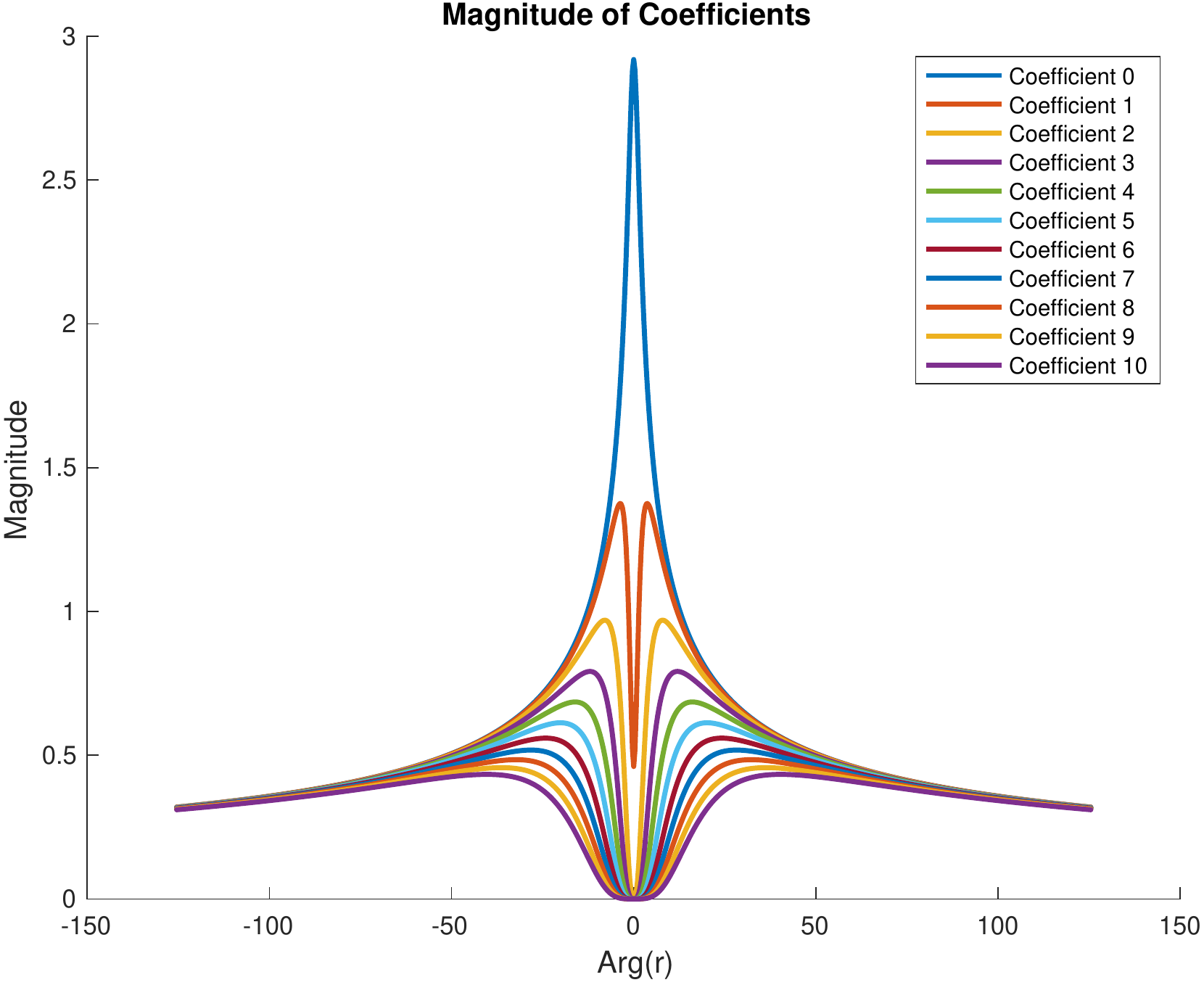}
\end{center}
\caption{A plot of the magnitudes of the first 11 coefficients for
$r=1+ir_y$ for $-20\pi \le r_y\le 20\pi$.}\label{fig:Magnitudes}
\end{figure}

Pragmatically, to find the values of $r$ for which the $a_0(r)$ are purely
real or imaginary requires solving Equation \eqref{eqn:Ak0Imag2}
numerically and then checking the conditions in Equations \eqref{eqn:Ak0Imagl}
and \eqref{eqn:Ak0Real}.

\subsection{The Magnitudes of the $a_k(r)$}

Numerical work by \cite{Pugh2004} showed that the magnitudes of the
Lanczos coefficients decreased rapidly for increasing $k$,
 motivating us to examine this property in the complex case. 
Figure (\ref{fig:Magnitudes}) presents the magnitudes of the first
11 Lanczos coefficients for $r=1+ir_y$ with $-20\pi \le y\le 20\pi$.
It is clear from this Figure that the coefficients with $r$ purely real
are turning points (maximum, minimum or local minimum) in the
trajectory of the $a_k(r)$ through complex space. We observed this
phenomena for all real values of $r$ we examined, 
$\mathcal{R}e(r)=1,3,5,7,9,11,13$.

\begin{figure}[h]
\begin{center}
\includegraphics[width=12cm]{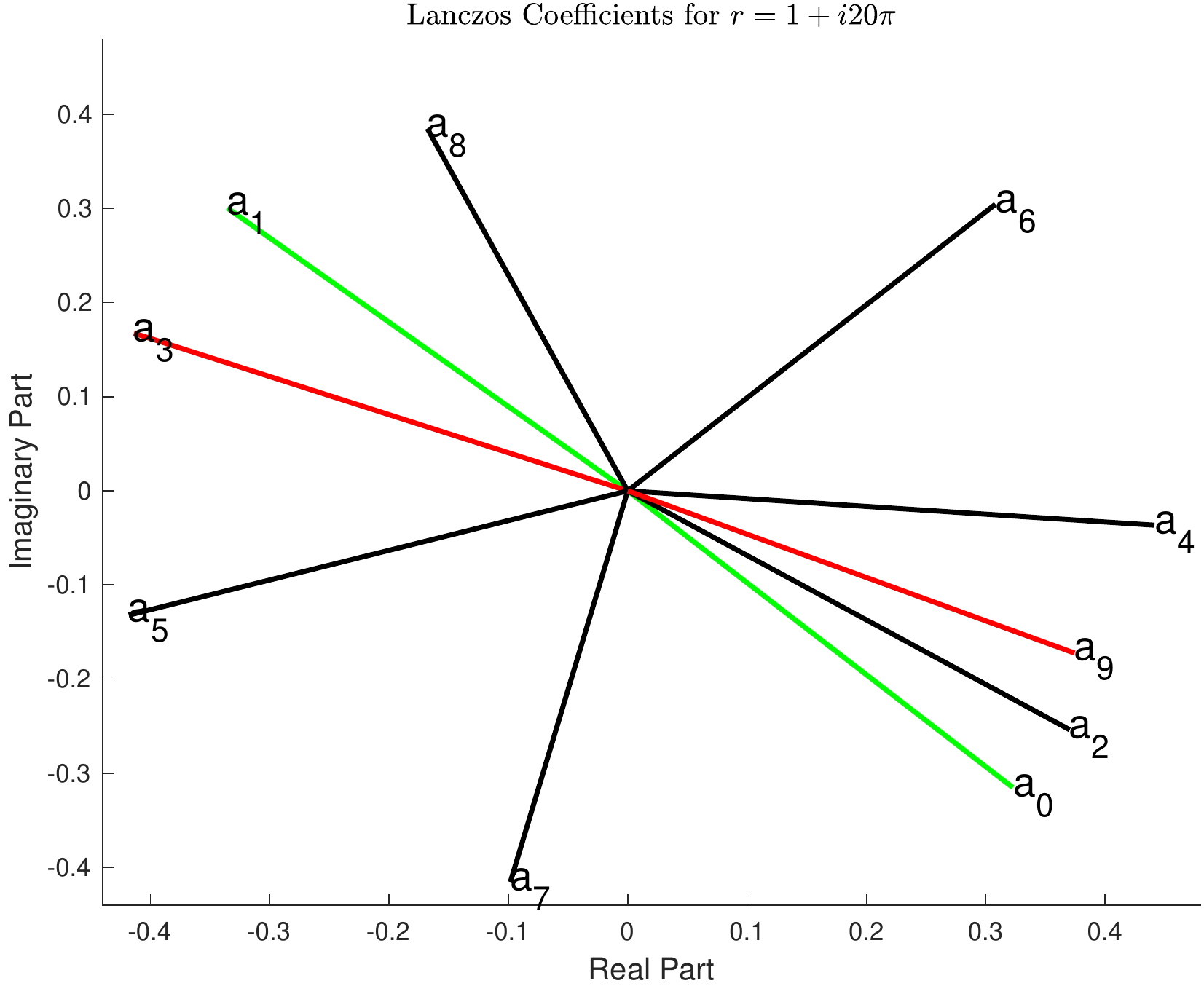}
\end{center}
\caption{A plot of the locations in the complex plane
of the first 10 Lanczos coefficients for
$r=1+20\pi i$.}\label{fig:R1p20Pi}
\end{figure}

Further, as $|\mathcal{I}m(r)|$ increases the magnitude of the coefficients
converge. This is illustrated in Figure (\ref{fig:R1p20Pi}) in which
we present the location in the complex
plane of the first ten Lanczos coefficients with $r=1+20\pi i$i. Their numerical
values are presented in  
Table (\ref{tab:R1Pi20}). We used these in one
set of numerical evaluations discussed in Section (\ref{sec:Numerical})
below. Both \cite{lanczos1964} and \cite{Pugh2004} noted that for real $r$
the signs of the Lanczos coefficients alternate with $a_k(r)>0$ for $k$
even and $a_k(r) <0$ for $k$ odd. Although we did not do extensive numerical
computations to large values of $k$, in the work we did do, we never found two
successive $a_k(r)$ in the same quadrant. 

\begin{table}
\begin{center}
\begin{tabular}{l|r}
\hline
Coefficient & Numerical Value \\
\hline
$a_0$  &  0.32272800 - 0.31511539i \\
$a_1$  & -0.33566576 + 0.30053412i \\
$a_2$  &  0.37045027 - 0.25375121i \\
$a_3$  & -0.41388748 + 0.16751781i \\
$a_4$  &  0.44150035 - 0.03647994i \\
$a_5$  & -0.41823596 - 0.13196358i \\
$a_6$  &  0.30817352 + 0.30446707i \\
$a_7$  & -0.09837046 - 0.41551337i \\
$a_8$  & -0.16793546 + 0.38494230i \\
$a_9$  &  0.37447472 - 0.172296311i \\
\hline
\end{tabular}
\end{center}
\caption{The numerical values of the first ten Lanczos coeffients
with free parameter $r=1+i20\pi$, corresponding to the those in
Figure (\ref{fig:R1p20Pi}). Their performance in the
approximation is  discussed briefly in Section 
(\ref{sec:Numerical}).}\label{tab:R1Pi20}
\end{table}

%% file: Numerical.tex
\section{Numerical Evaluation}\label{sec:Numerical}

It is obvious that using a Lanczos approximation with
complex $r$ increases the computational
cost of calculating the $\Gamma$-function when the same number of 
coefficients are used. Thus to use complex coefficients
requires a compelling case that there are significant advantages
in the convergence rate and/or accuracy of the approximation
compared to real coefficients. In this section we examine whether
such a case can be made.

\begin{figure}[h]
\begin{center}
\includegraphics[width=12cm]{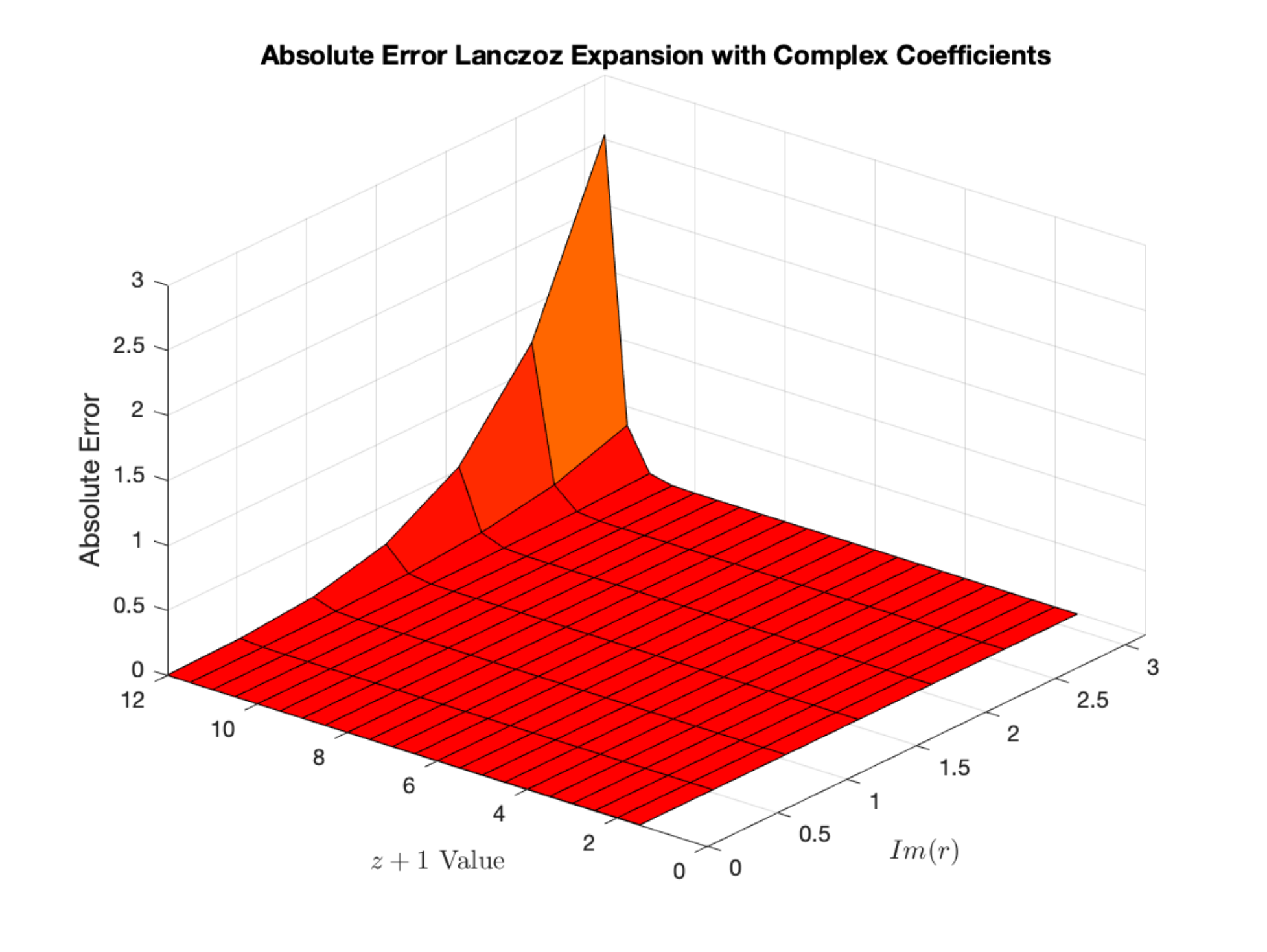}
\end{center}
\caption{The absolute error of the Lanczos $\Gamma$ approximation for $z+1$ 
between 1.5 and 12 using 10 terms with $r=1+ir_y$ with $0\le r_y \le 2\pi$.
}\label{fig:ComplexRAbs}
\end{figure}

\begin{figure}[h]
\begin{center}
\includegraphics[width=12cm]{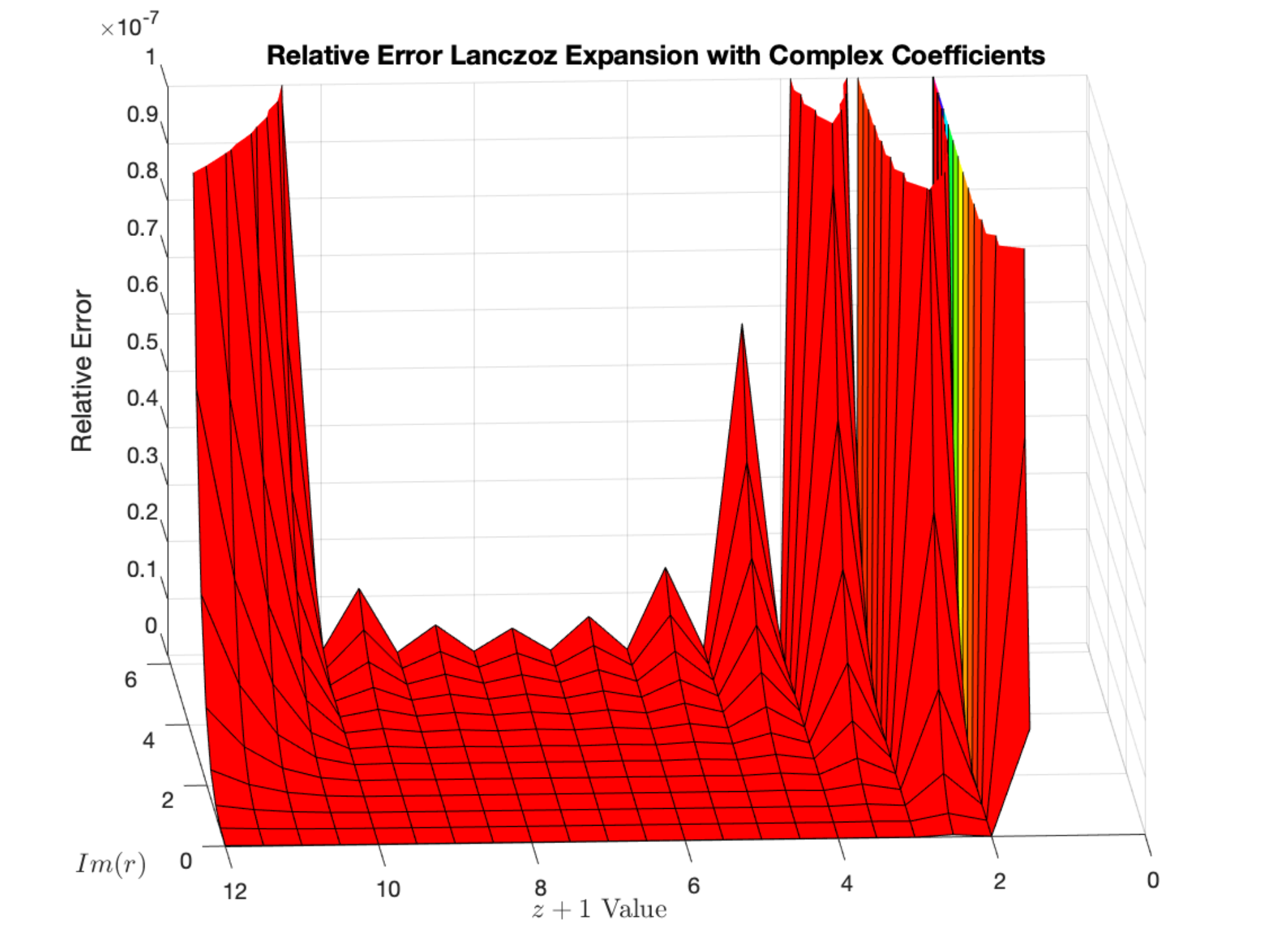}
\end{center}
\caption{The relative error of the Lanczos $\Gamma$ approximation for $z+1$
between 1.5 and 12 using 10 terms with $r=1+ir_y$ with $0\le r_y \le 2\pi$. The
top of the graph is truncated at an error of one part in
$10^7$.}
\label{fig:ComplexRRel}
\end{figure}

\begin{figure}[h]
\begin{center}
\includegraphics[width=12cm]{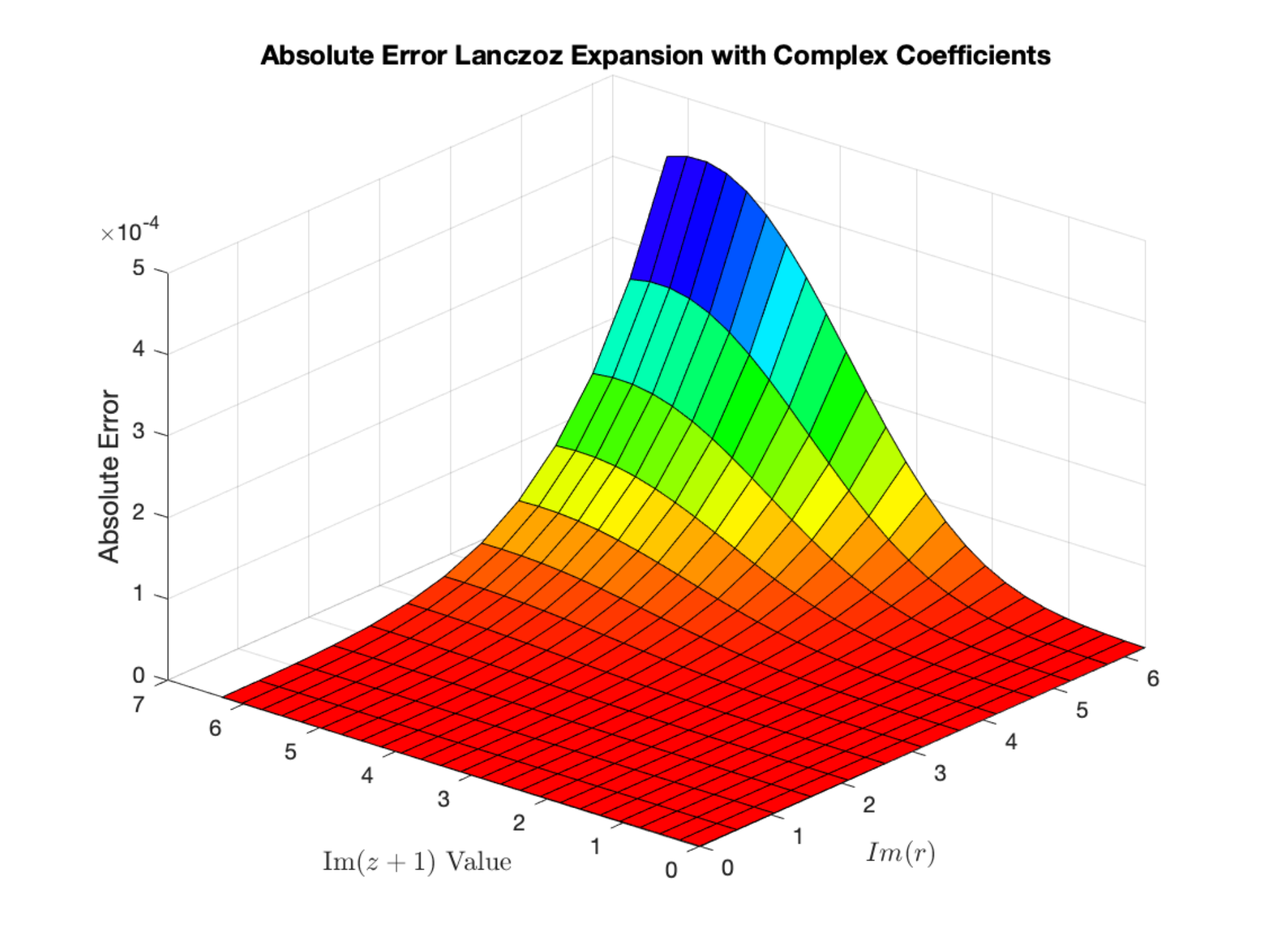}
\end{center}
\caption{The absolute error of the Lanczos $\Gamma$ approximation for 
$z+1=6+iy$ with $0\le y \le 2\pi$
using 10 terms with $r=1+ir_y$ with $0\le r_y \le 2\pi$.
}\label{fig:ComplexZAbs}
\end{figure}

\begin{figure}[h]
\begin{center}
\includegraphics[width=12cm]{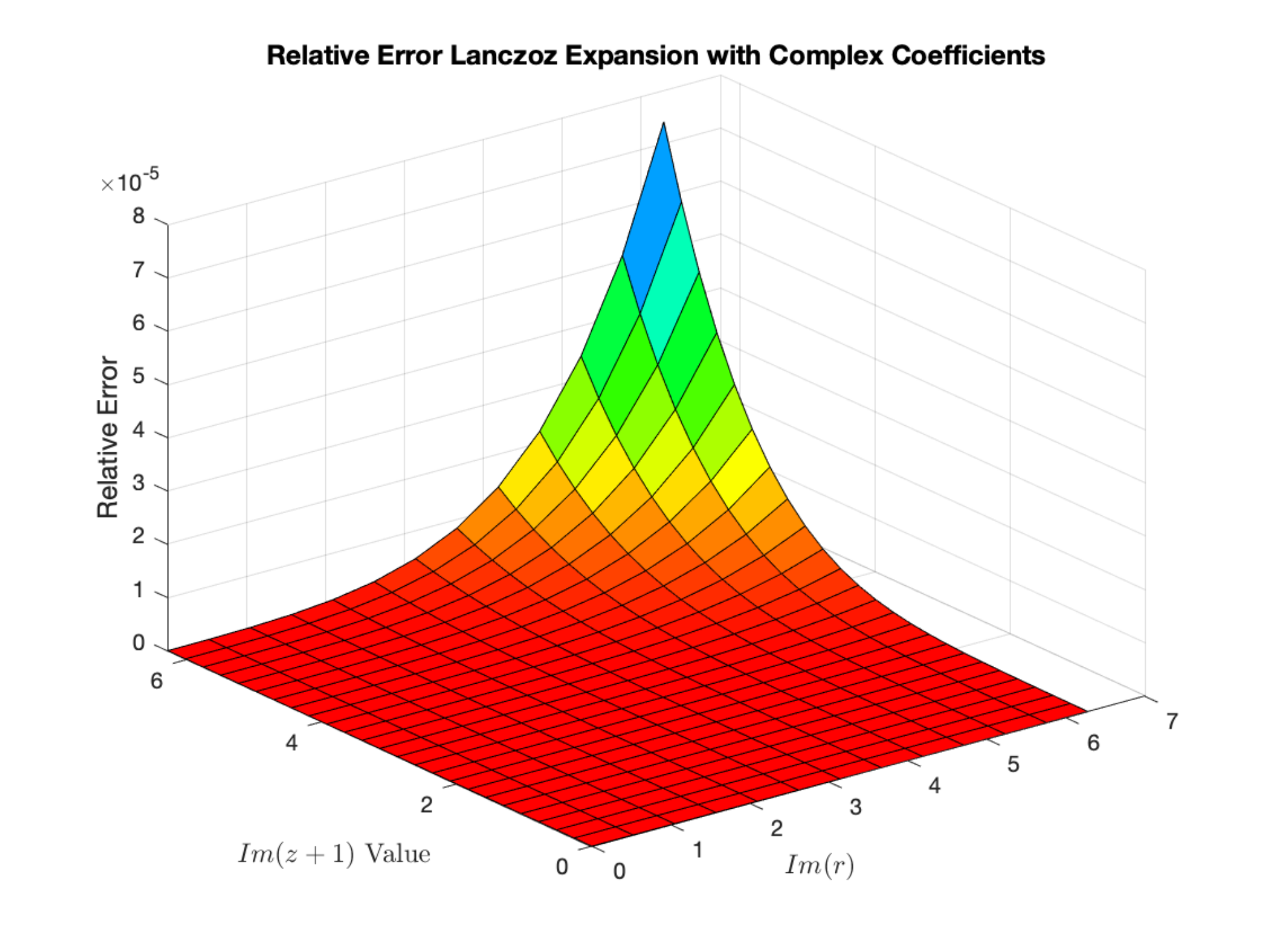}
\end{center}
\caption{The relative error of the Lanczos $\Gamma$ approximation for 
$z+1=6+iy$ with $0\le y \le 2\pi$
using 10 terms with $r=1+ir_y$ with $0\le r_y \le 2\pi$.}
\label{fig:ComplexZRel}
\end{figure}

In all cases we used 10 coefficients in the approximation.
Figures (\ref{fig:ComplexRAbs}) through (\ref{fig:ComplexZRel}) 
present the results, we discuss how the data in each Figure 
was generated and comment on its significance.

We generated 13 sets of 10 coefficients using the recursive method 
described above by setting $r=1+in\pi/6$
where $n=0,1,\ldots,12$. In Figures (\ref{fig:ComplexRAbs}) 
and (\ref{fig:ComplexRRel}) the imaginary
part of $r$ is presented on the axis label $Im(r)$. We generated 22
values of $z+1=1.5,2.0,2.5,\ldots,12$ for which exact values could
be calculated. These were obtained for integer values
from the factorial property, namely
$\Gamma(n+1)=n!$ and for half integer values from the functional
equation 
$$
\Gamma(z+1)=z\Gamma(z)
$$
coupled with the value $\Gamma(1/2)=\sqrt{\pi}$.

Figure (\ref{fig:ComplexRAbs}) presents the absolute value of the
error while Figure (\ref{fig:ComplexRRel}) presents the relative error,
that is the absoluate error is divided by the exact value. In both
cases it is clear that the Lanczos coefficients generated by
real $r$ performed better. The absolute error is never large (less
than 3 for $11! \approx 3.9\times 10^7$ but largest error occurs for
coefficients generated with  $r=1+i2\pi$.

The relative error presented in Figure (\ref{fig:ComplexRRel}) is more
informative. The vertical axis is truncated at an error of $10^{-7}$ in
order to better see the relative error at the half-integer values..
It clear that in all cases the relative error is indistiquishable from
zero for integer values up to $z+1=10$ for all values of$r$. 
This is because they were used to
generate the coefficients. However, the half integer values of $z+1$ 
show increasing relative errors as the imaginary part increases from
$0$ to $2\pi$. Somewhat similar, once we progress beyond the integer
values of $z+1$ which we used to generate the coefficients, i.e. $z+1=10$,
the relative error builds rapidly with increasing imaginary part of $r$.
At least on the scale used in the Figure the relative error for the
coefficients are indistinquiable from zero for real $r$.

For Figures (\ref{fig:ComplexZAbs}) and (\ref{fig:ComplexZRel}), using
the same coefficients as above,  we considered
the performance of the approxmation for complex values of $z+1$. We 
set $z+1=6+n\pi i/12$ for $n=0,1,\ldots,24$. To find  comparison values
we used the \cite{Godfrey2001} implementation discussed above. 
For both the absolute and relative error cases, the Lanczos
approximation performance
deteriorated as the imaginary part of $z+1$ increased and as the
imaginary part of $r$ increased together. 

Figure (\ref{fig:Magnitudes}) suggests that as the imagary part of $r$ increases
all coefficients converge to the same magnitude. For the coefficients 
displayed in Figure (\ref{fig:R1p20Pi}), i.e.\ for $r=1+20\pi i$, we
tested its performance against the \cite{Godfrey2001} implementation. We do
not present detailed results here, suffice to say that as the imaginary
part of $z+1$ increased the relative error increased very rapidly.

%% file: Conclusions.tex
\section{Conclusions}\label{sec:Conclude}

We explored the possibility that complex values of the free parameter $r$ 
in the
\cite{lanczos1964} approximation for the $\Gamma$-function
could improve the rate of convergence of the approximation. However, we
must conclude that there is no advantage to using complex values of
$r$, in fact, it appears uniformly detrimental. The real values of $r$
for the values we explored uniformly gave the lowest relative error
when compared with exact values of the $\Gamma$-function at integer
and half integer values. Further, Lanczos coefficients generated
by real values of $r$ also uniformly gave the lowest relative error
for complex values of the $z$ in $\Gamma(z+1)$. 
 
Unless there are purely mathematical reasons for further
studying the properties
of complex Lanczos coefficients they appears to be a dead end in applied
numerical analysis.

%% file: AllGraphs.tex
\section{Plots of Coefficients 0 through 10}\label{sec:Supplement}

\subsection{Real and Imaginary Parts of the Coefficients}

\begin{figure}[h]
\begin{center}
\includegraphics[width=10cm]{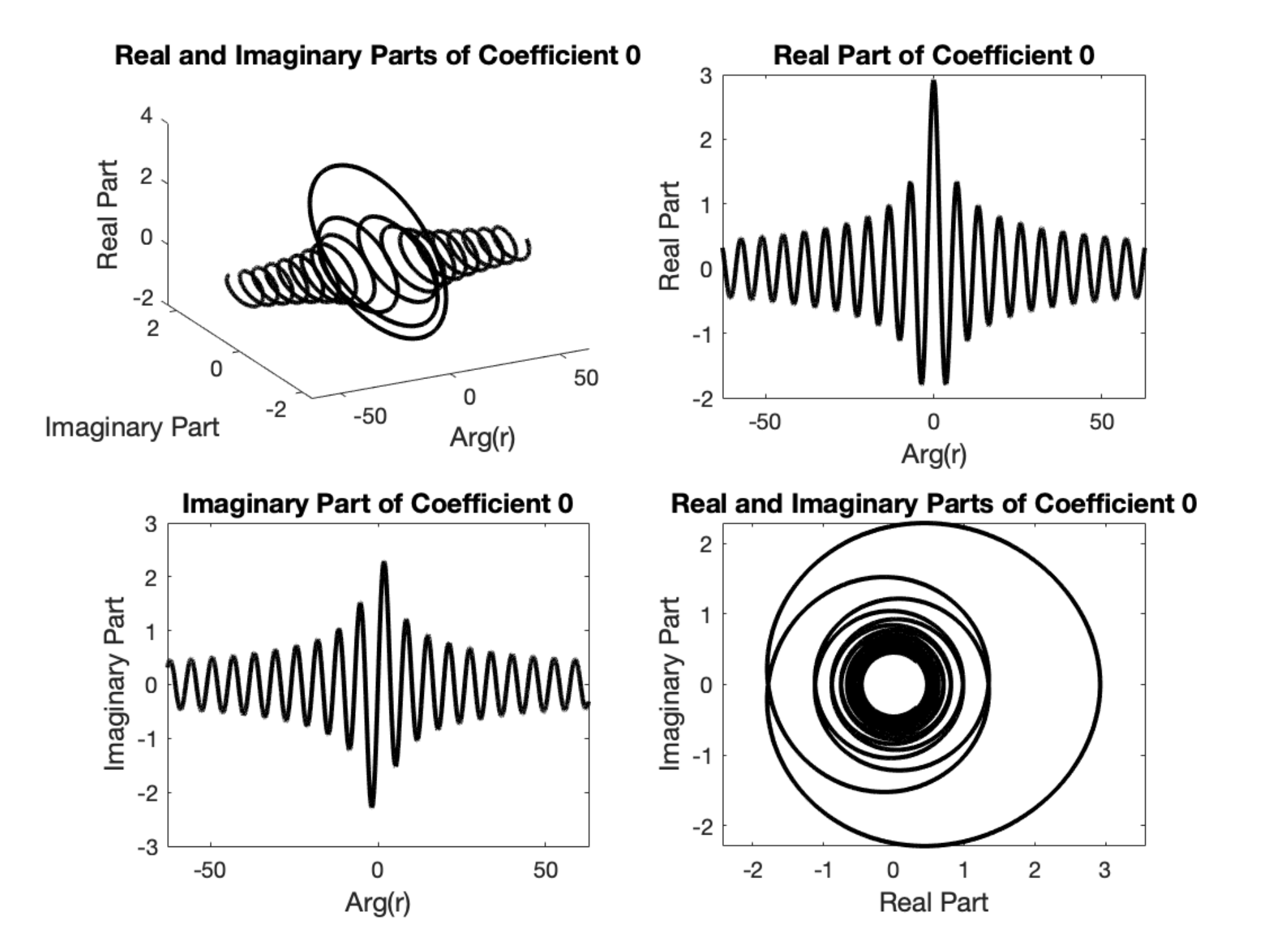}
\end{center}
\caption{Four displays of the real and imaginary parts for coefficient
$a_0(r)$ for $r=1+ir_y$ for $-20\pi \le r_y\le 20\pi$.}\label{fig:A04Plot}
\end{figure}

\begin{figure}
\begin{center}
\includegraphics[width=10cm]{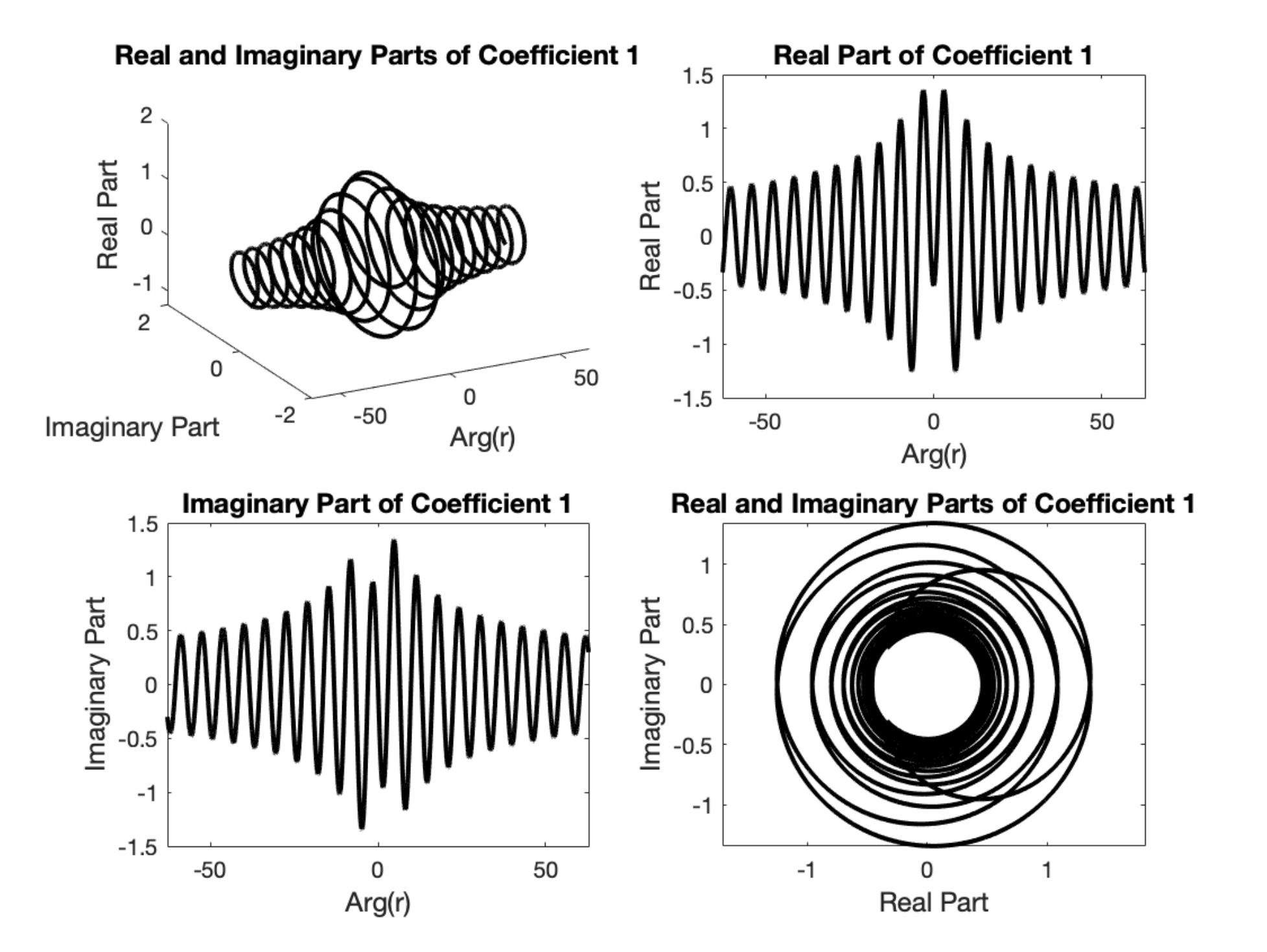}
\end{center}
\caption{Four displays of the real and imaginary parts for coefficient
$a_1(r)$ for $r=1+ir_y$ with  $-20\pi \le r_y\le 20\pi$.}\label{fig:A14Plot}
\end{figure}

\begin{figure}
\begin{center}
\includegraphics[width=10cm]{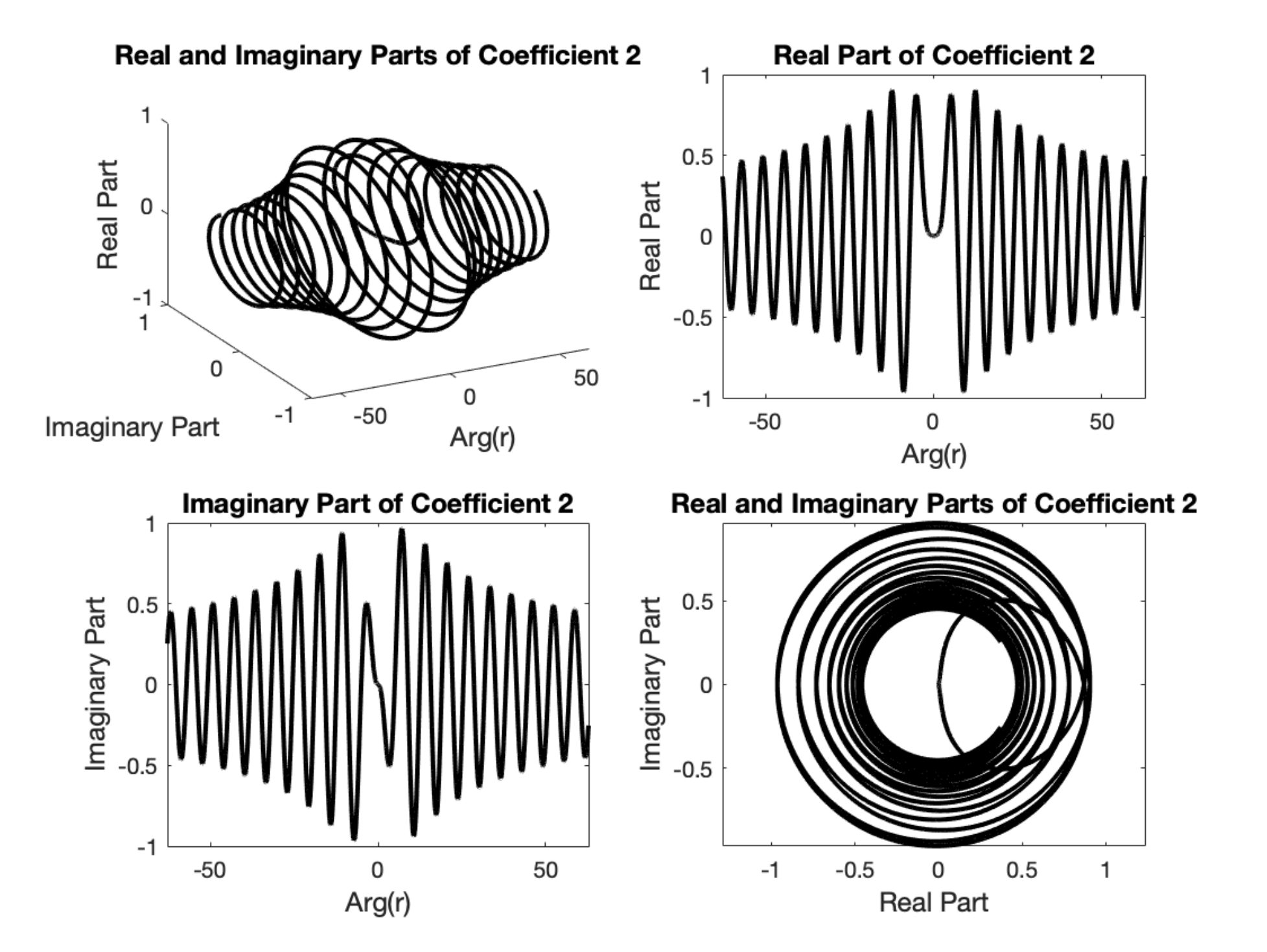}
\end{center}
\caption{Four displays of the real and imaginary parts for coefficient
$a_2(r)$ for $r=1+ir_y$ with  $-20\pi \le r_y\le 20\pi$.}\label{fig:A24Plot}
\end{figure}

\begin{figure}
\begin{center}
\includegraphics[width=10cm]{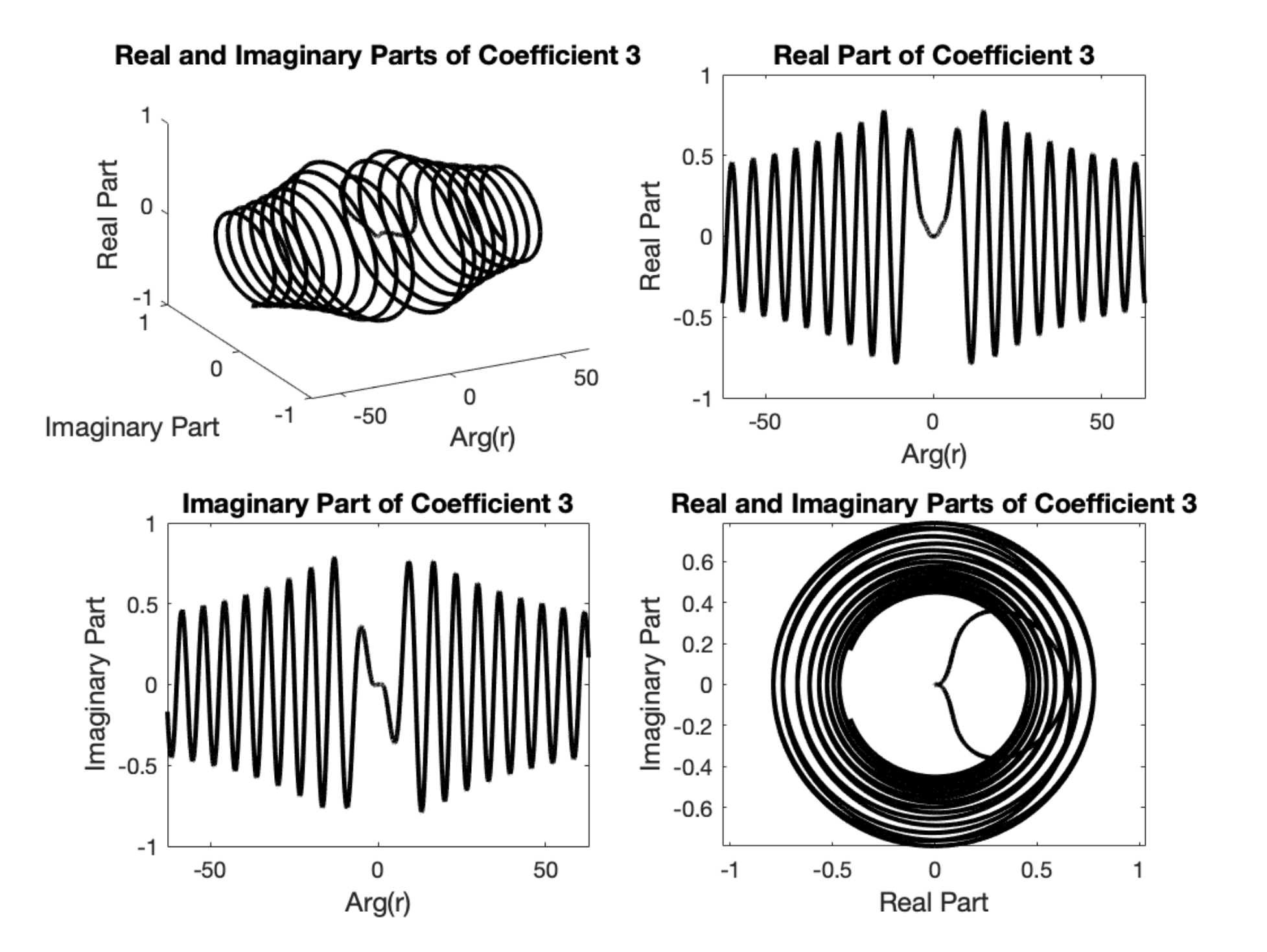}
\end{center}
\caption{Four displays of the real and imaginary parts for coefficient
$a_3(r)$ for $r=1+ir_y$ with  $-20\pi \le r_y\le 20\pi$.}\label{fig:A34Plot}
\end{figure}

\begin{figure}
\begin{center}
\includegraphics[width=10cm]{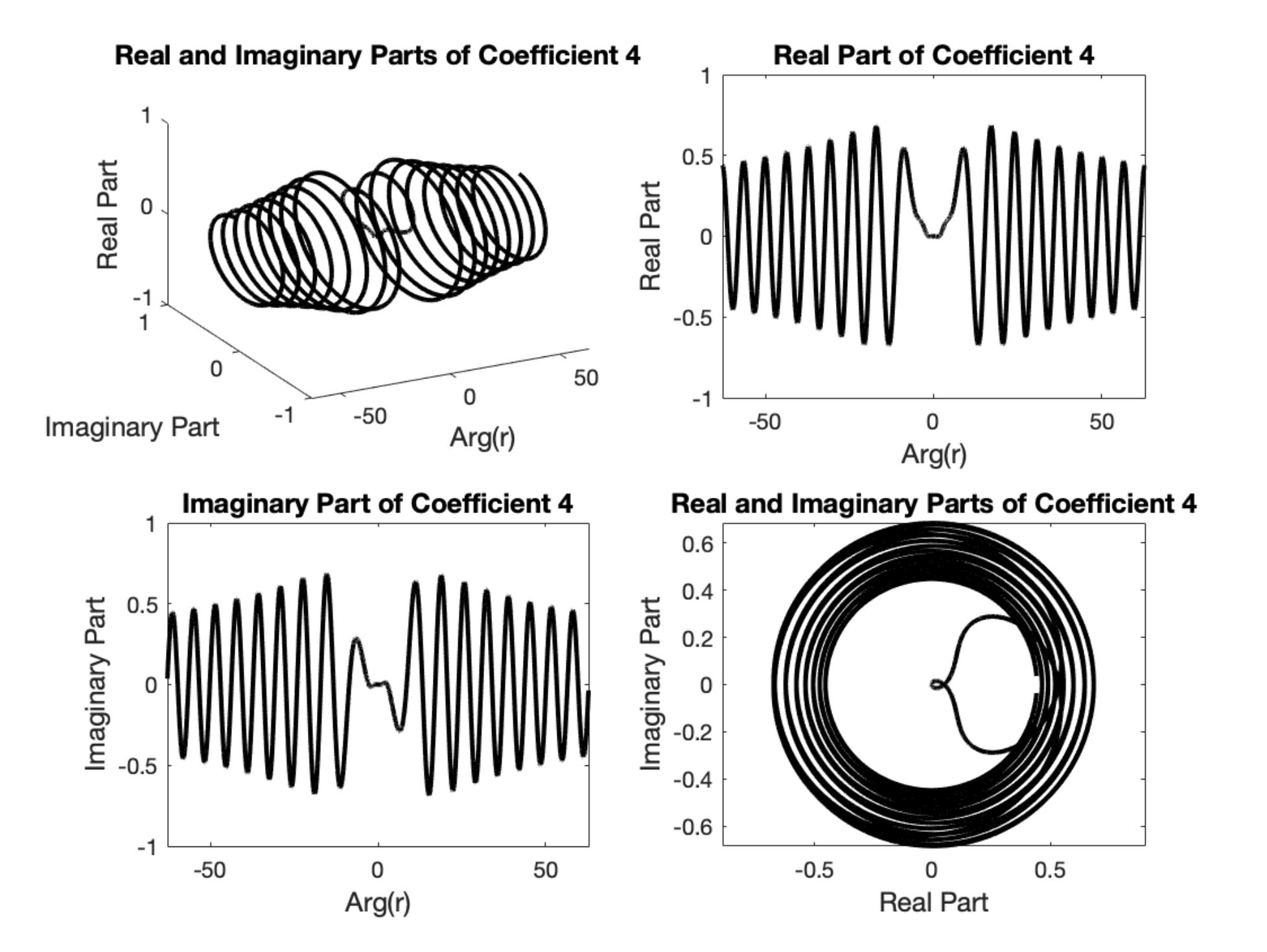}
\end{center}
\caption{Four displays of the real and imaginary parts for coefficient
$a_4(r)$ for $r=1+ir_y$ with $-20\pi \le r_y\le 20\pi$.}\label{fig:A44Plot}
\end{figure}

\begin{figure}
\begin{center}
\includegraphics[width=10cm]{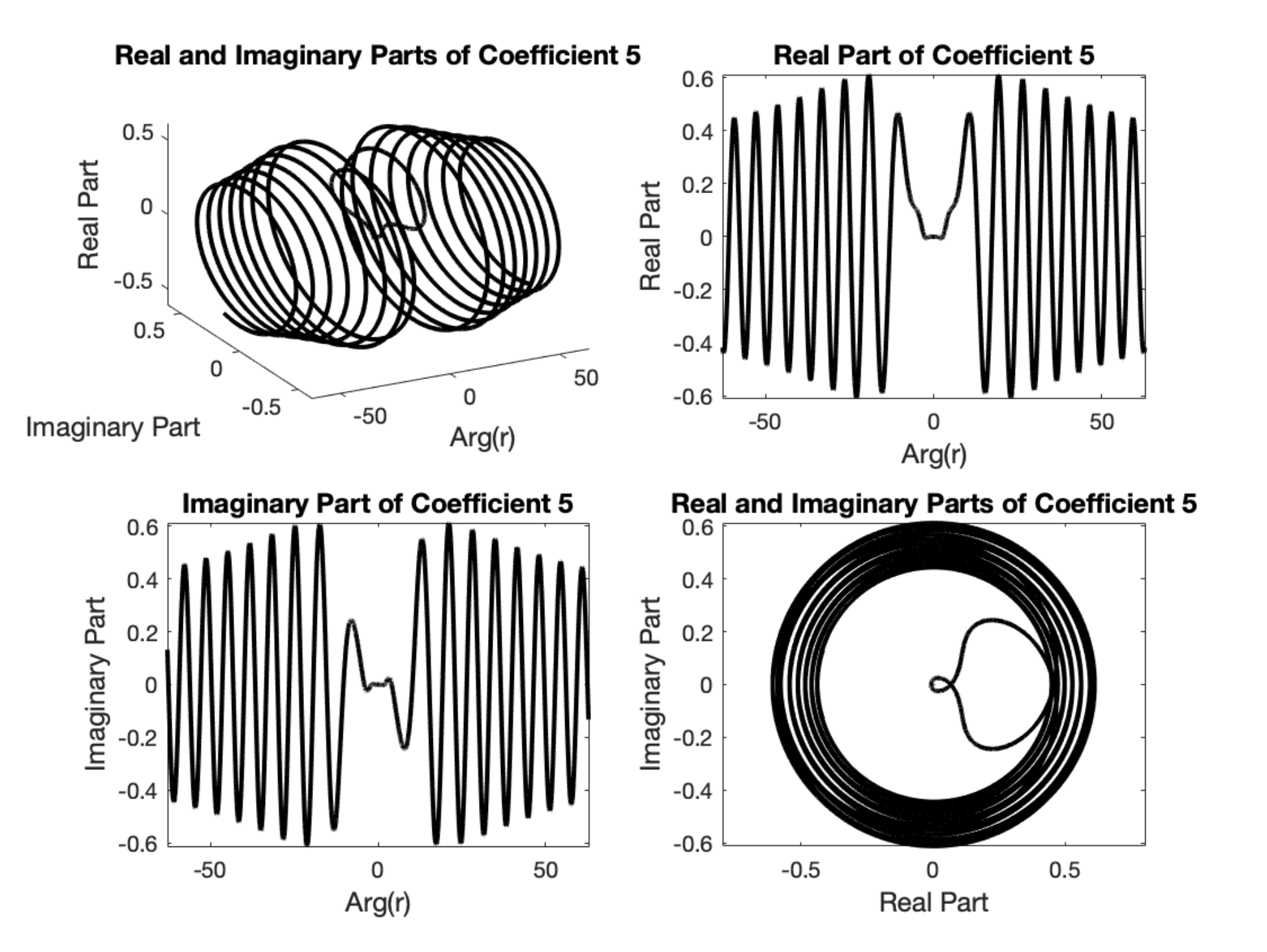}
\end{center}
\caption{Four displays of the real and imaginary parts for coefficient
$a_5(r)$ for $r=1+ir_y$ with $-20\pi \le r_y\le 20\pi$.}\label{fig:A54Plot}
\end{figure}

\begin{figure}
\begin{center}
\includegraphics[width=10cm]{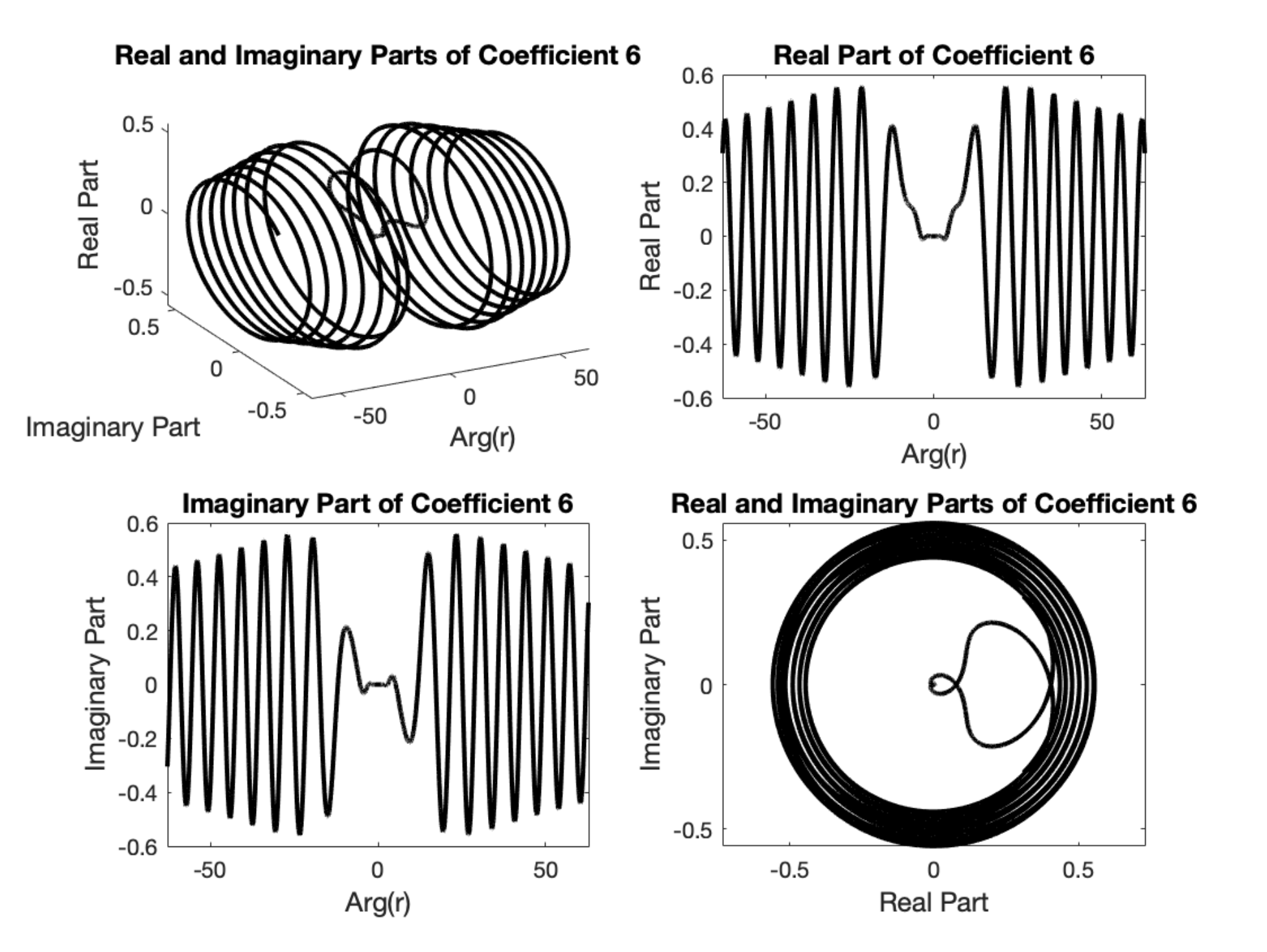}
\end{center}
\caption{Four displays of the real and imaginary parts for coefficient
$a_6(r)$ for $r=1+ir_y$ with $-20\pi \le r_y\le 20\pi$.}\label{fig:A64Plot}
\end{figure}

\begin{figure}
\begin{center}
\includegraphics[width=10cm]{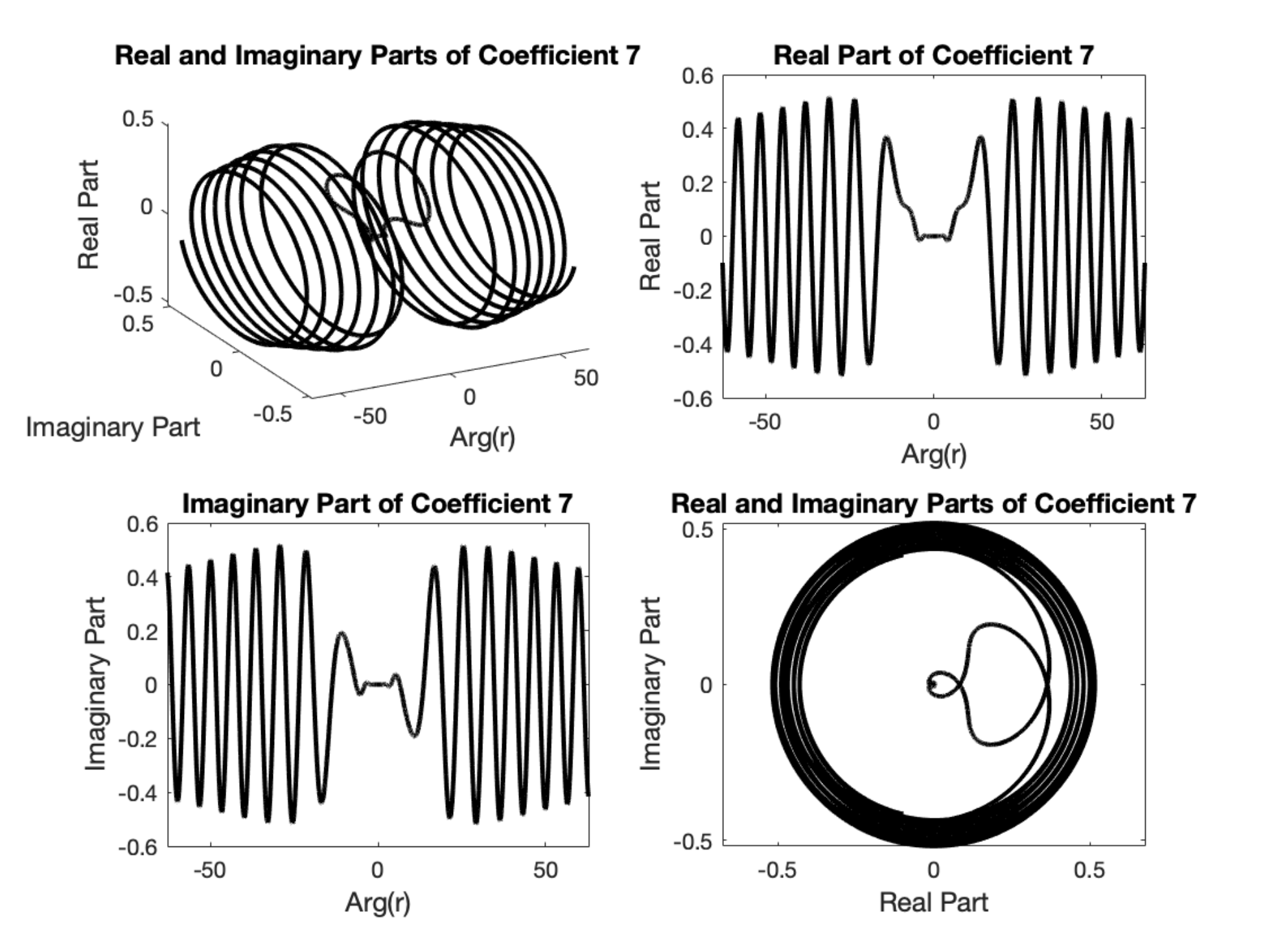}
\end{center}
\caption{Four displays of the real and imaginary parts for coefficient
$a_7(r)$ for $r=1+ir_y$ with $-20\pi \le r_y\le 20\pi$.}\label{fig:A74Plot}
\end{figure}

\begin{figure}
\begin{center}
\includegraphics[width=10cm]{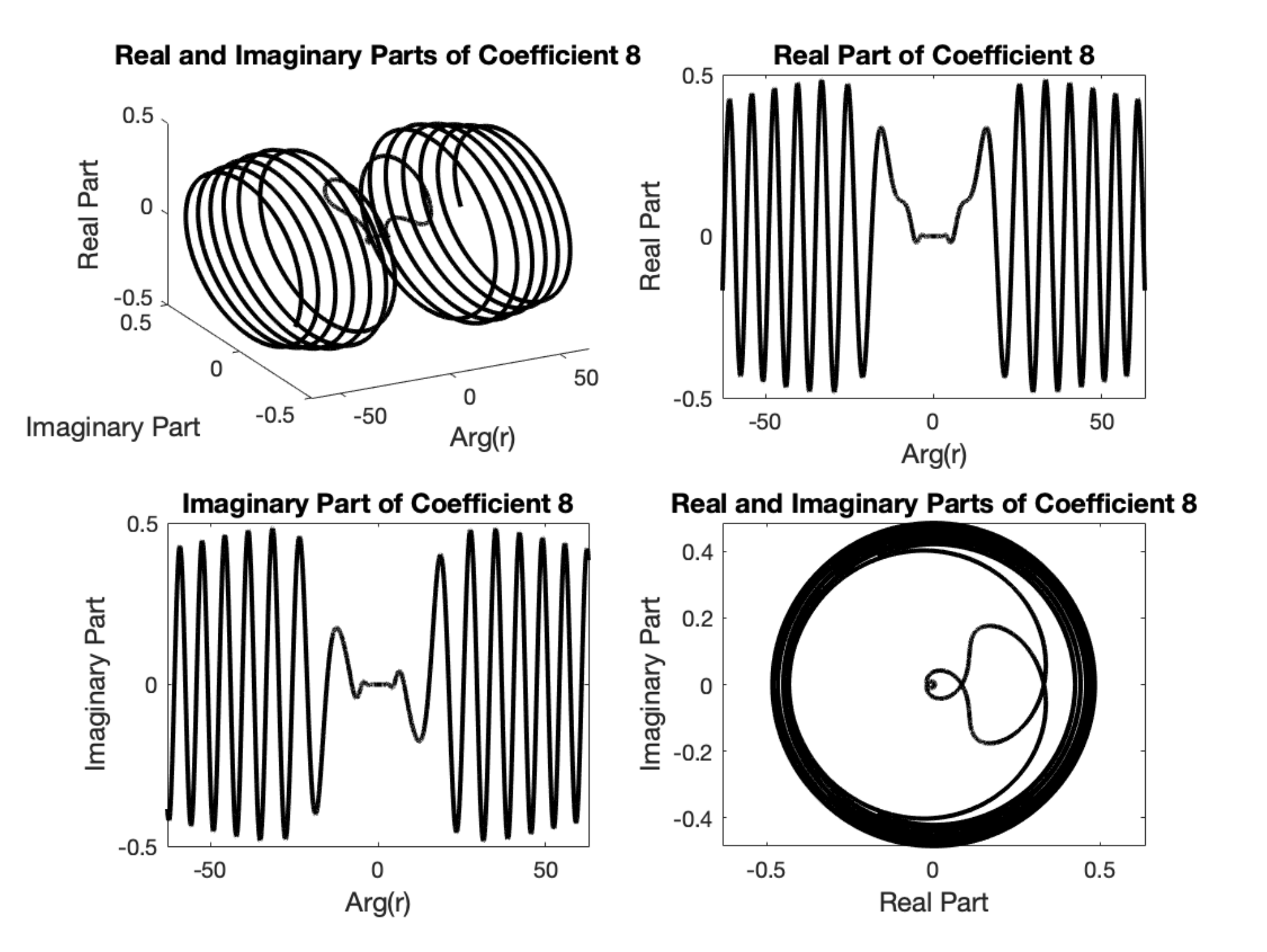}
\end{center}
\caption{Four displays of the real and imaginary parts for coefficient
$a_8(r)$ for $r=1+ir_y$ with $-20\pi \le r_y\le 20\pi$.}\label{fig:A84Plot}
\end{figure}

\begin{figure}
\begin{center}
\includegraphics[width=10cm]{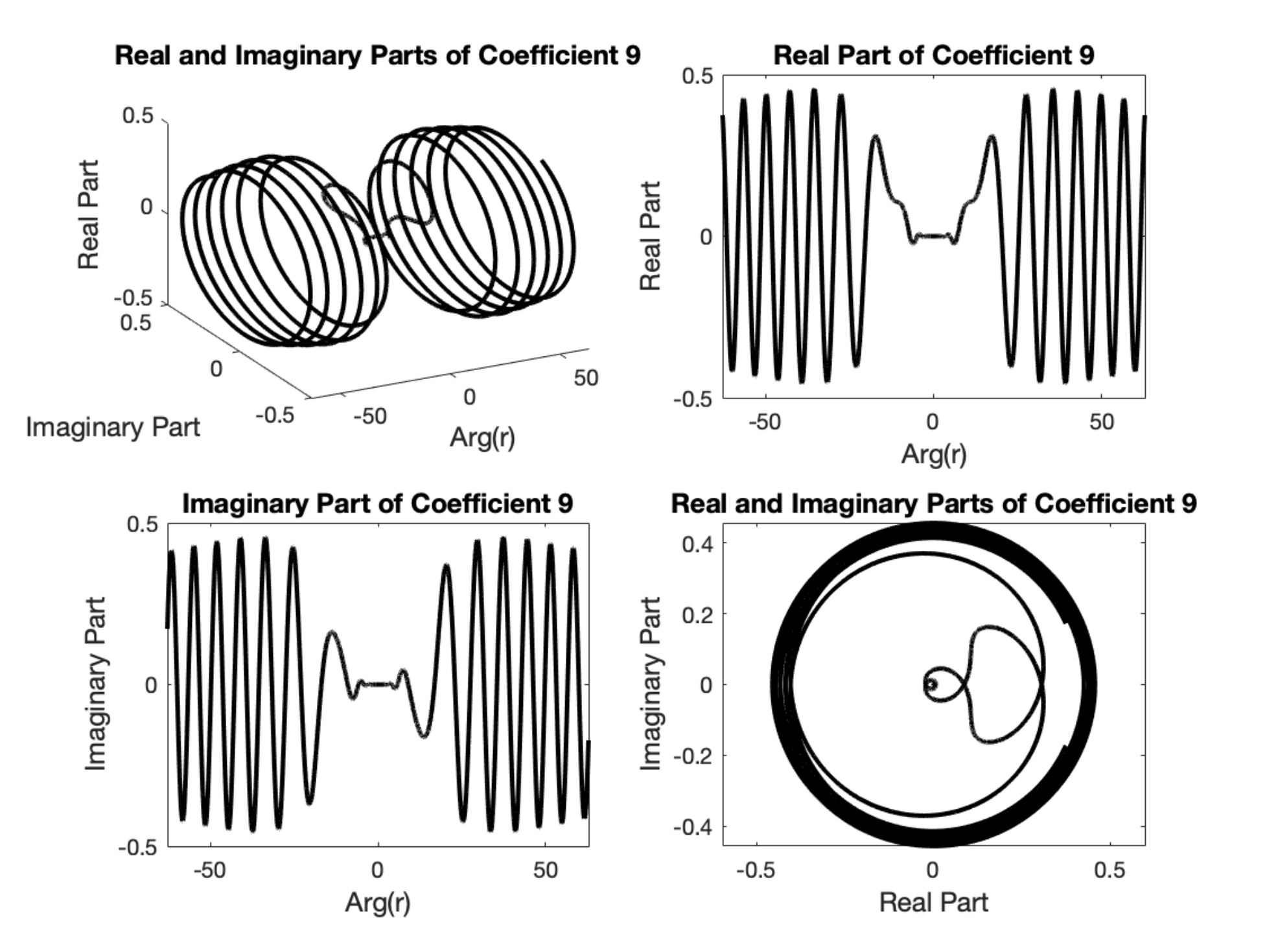}
\end{center}
\caption{Four displays of the real and imaginary parts for coefficient
$a_9(r)$ for $r=1+ir_y$ with $-20\pi \le r_y\le 20\pi$.}\label{fig:A94Plot}
\end{figure}

\begin{figure}
\begin{center}
\includegraphics[width=10cm]{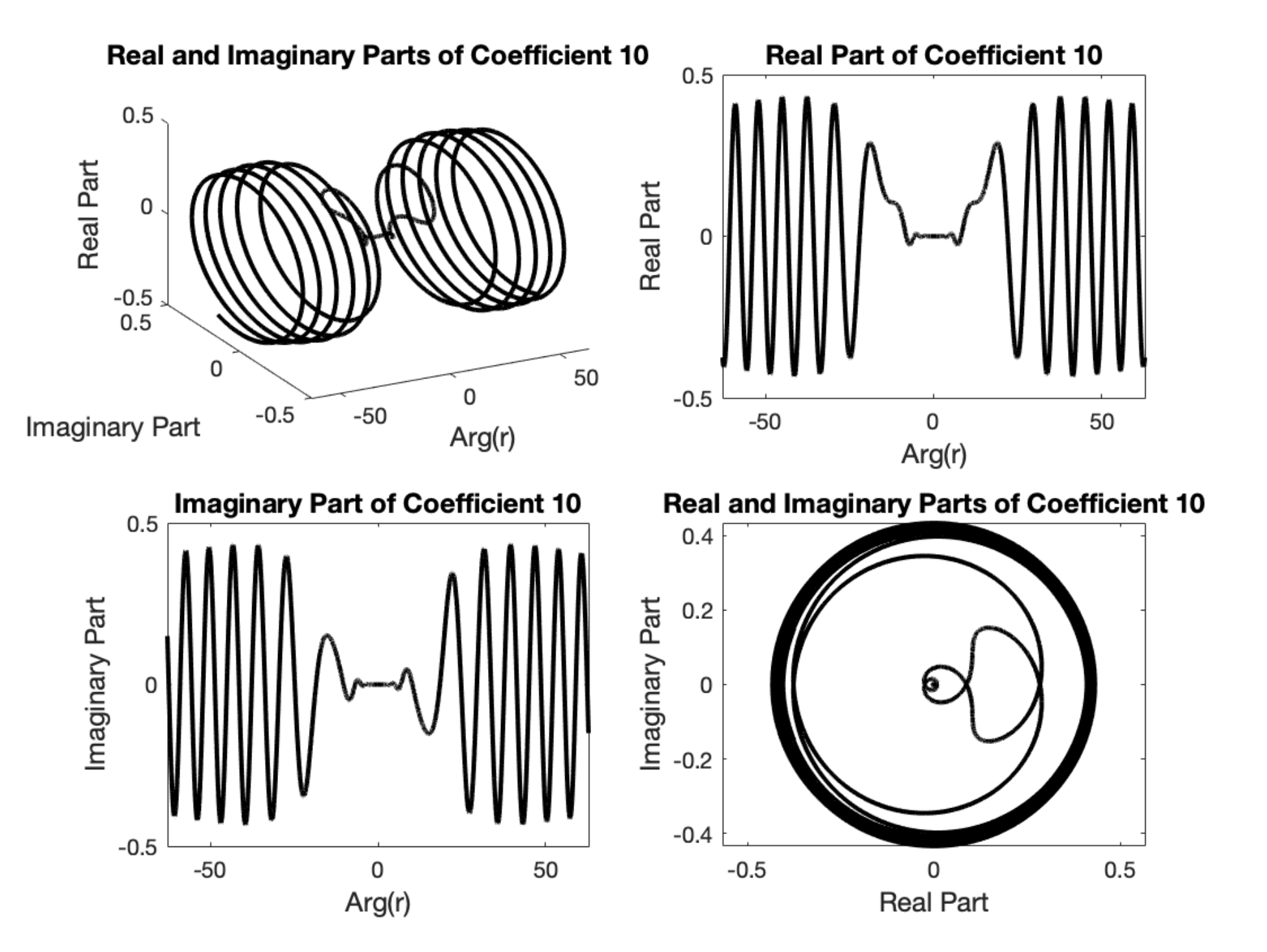}
\end{center}
\caption{Four displays of the real and imaginary parts for coefficient
$a_{10}(r)$ for $r=1+ir_y$ with $-20\pi \le r_y\le 20\pi$.}\label{fig:A104Plot}
\end{figure}

\clearpage

\subsection{Coefficients on the Riemann Sphere}


\begin{figure}[h]
\begin{center}
\includegraphics[width=10cm]{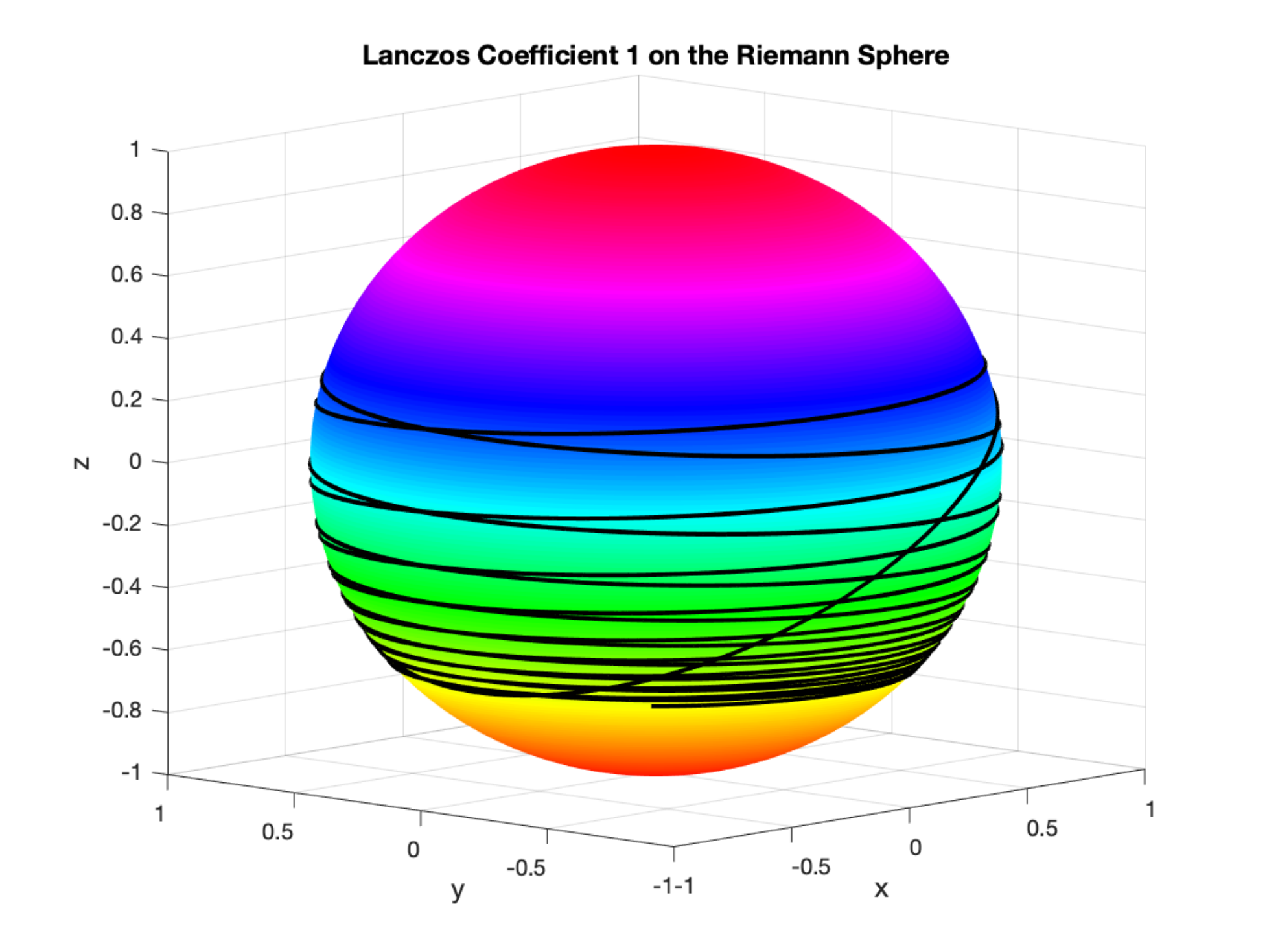}
\end{center}
\caption{The coefficient
$a_1(r)$ for $r=1+ir_y$ with  $-20\pi \le r_y\le 20\pi$ plotted on the
Riemann sphere.}\label{fig:A1RS}
\end{figure}

\begin{figure}
\begin{center}
\includegraphics[width=10cm]{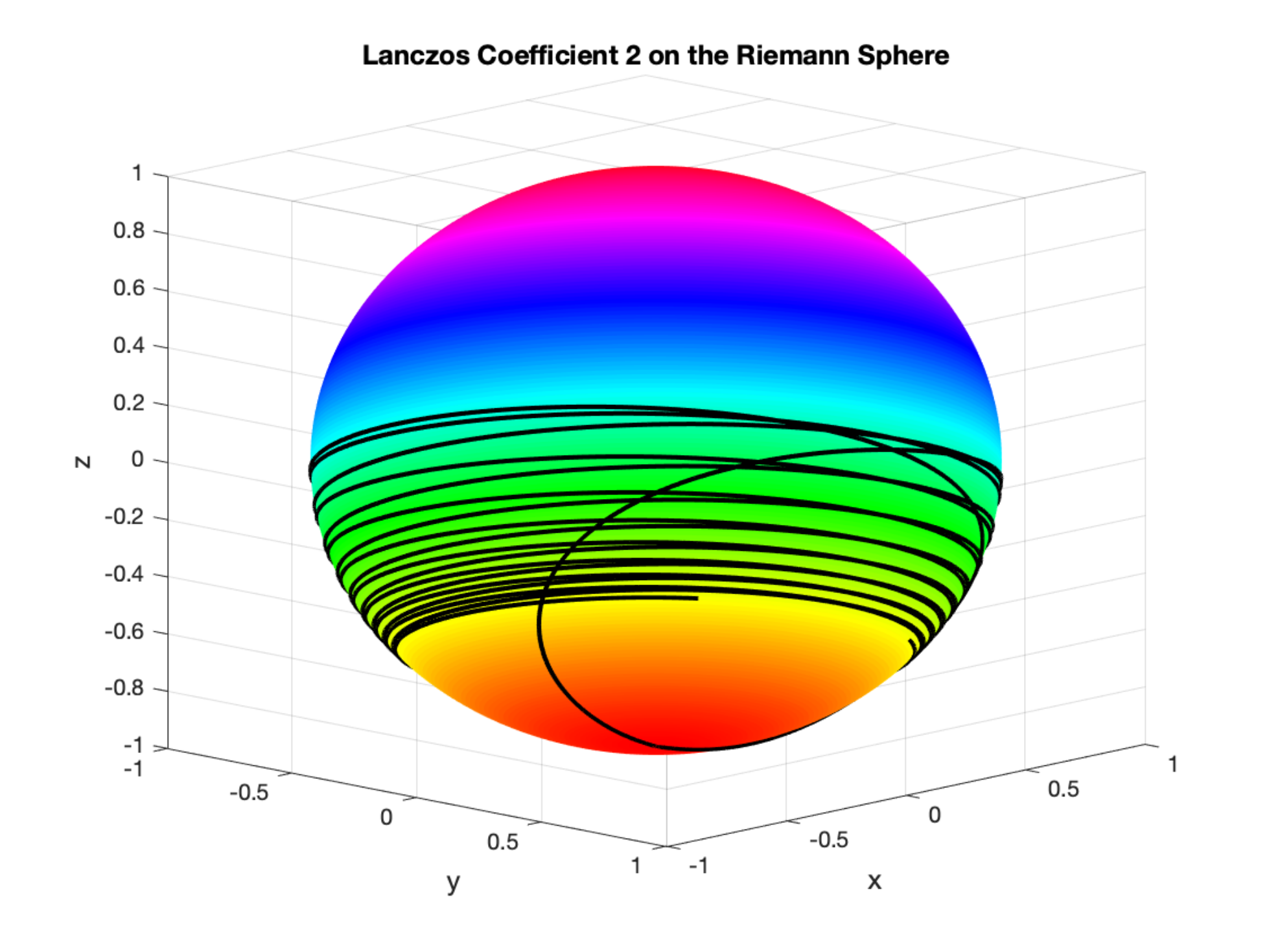}
\end{center}
\caption{The coefficient
$a_2(r)$ for $r=1+ir_y$ with  $-20\pi \le r_y\le 20\pi$ plotted on the
Riemann sphere.}\label{fig:A2RS}
\end{figure}

\begin{figure}
\begin{center}
\includegraphics[width=10cm]{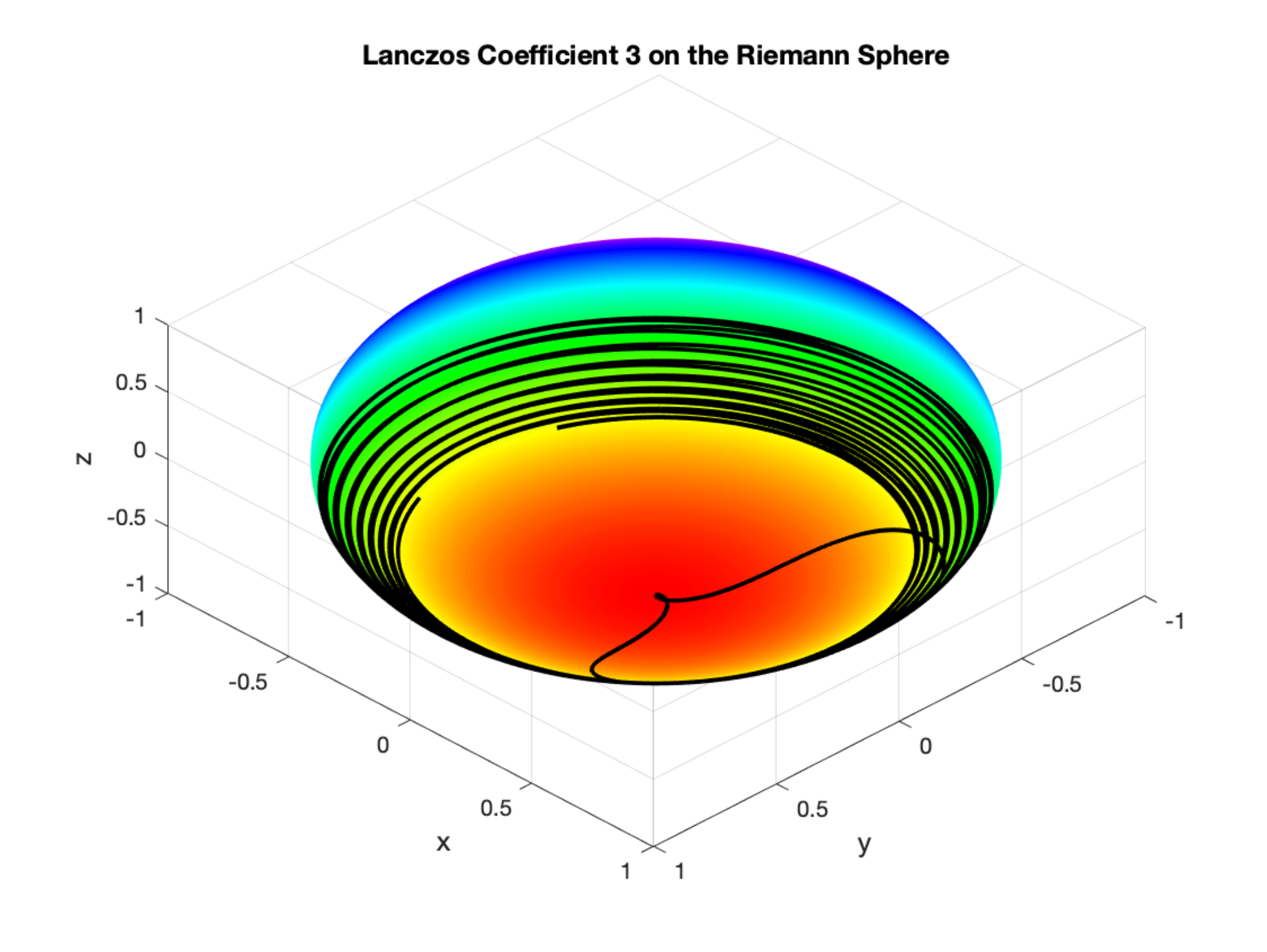}
\end{center}
\caption{The coefficient
$a_3(r)$ for $r=1+ir_y$ with  $-20\pi \le r_y\le 20\pi$ plotted on
the Riemann sphere.}\label{fig:A3RS}
\end{figure}

\begin{figure}
\begin{center}
\includegraphics[width=10cm]{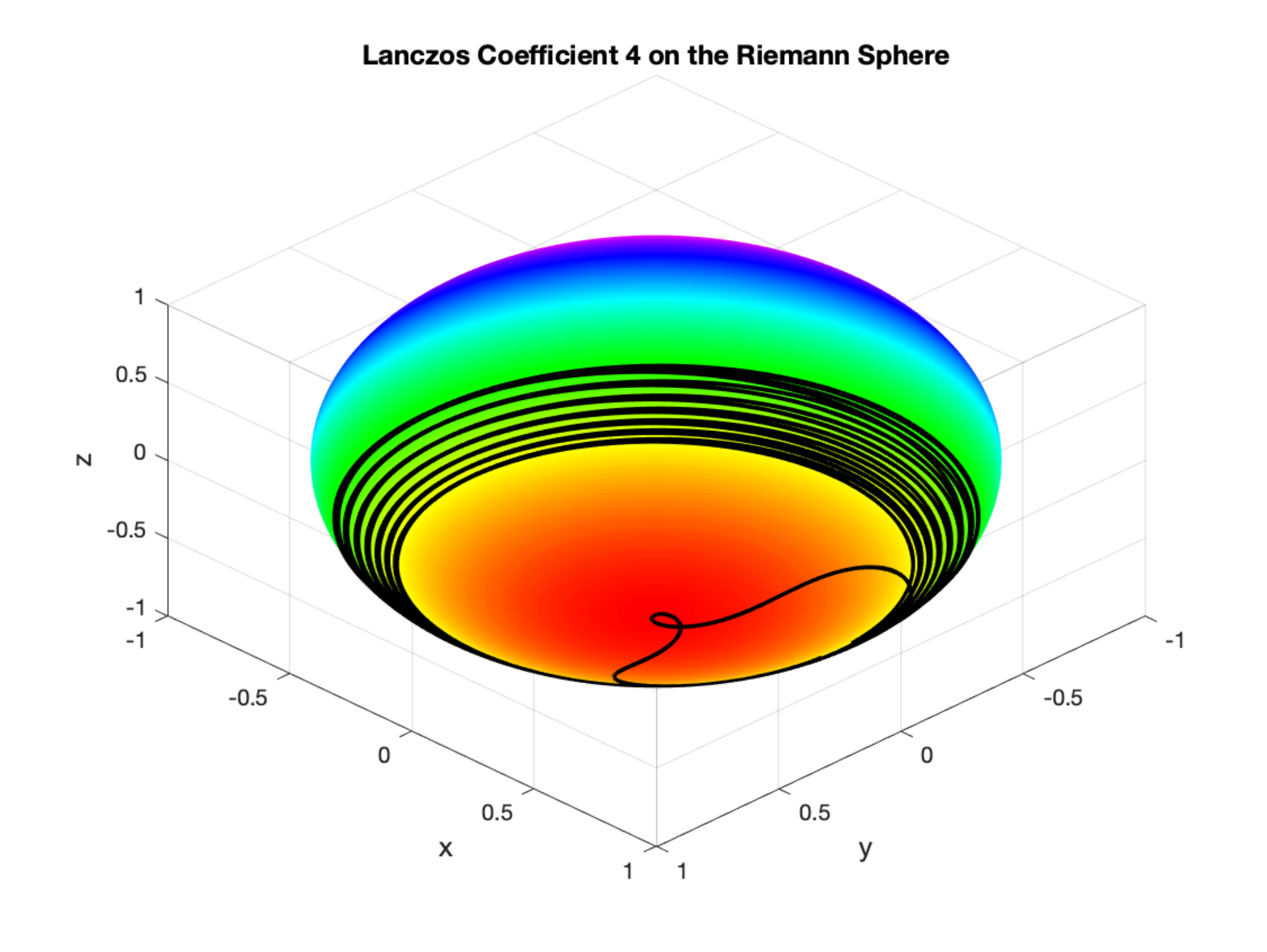}
\end{center}
\caption{The coefficient
$a_4(r)$ for $r=1+ir_y$ with $-20\pi \le r_y\le 20\pi$ plotted on the
Riemann sphere.}\label{fig:A4RS}
\end{figure}

\begin{figure}
\begin{center}
\includegraphics[width=10cm]{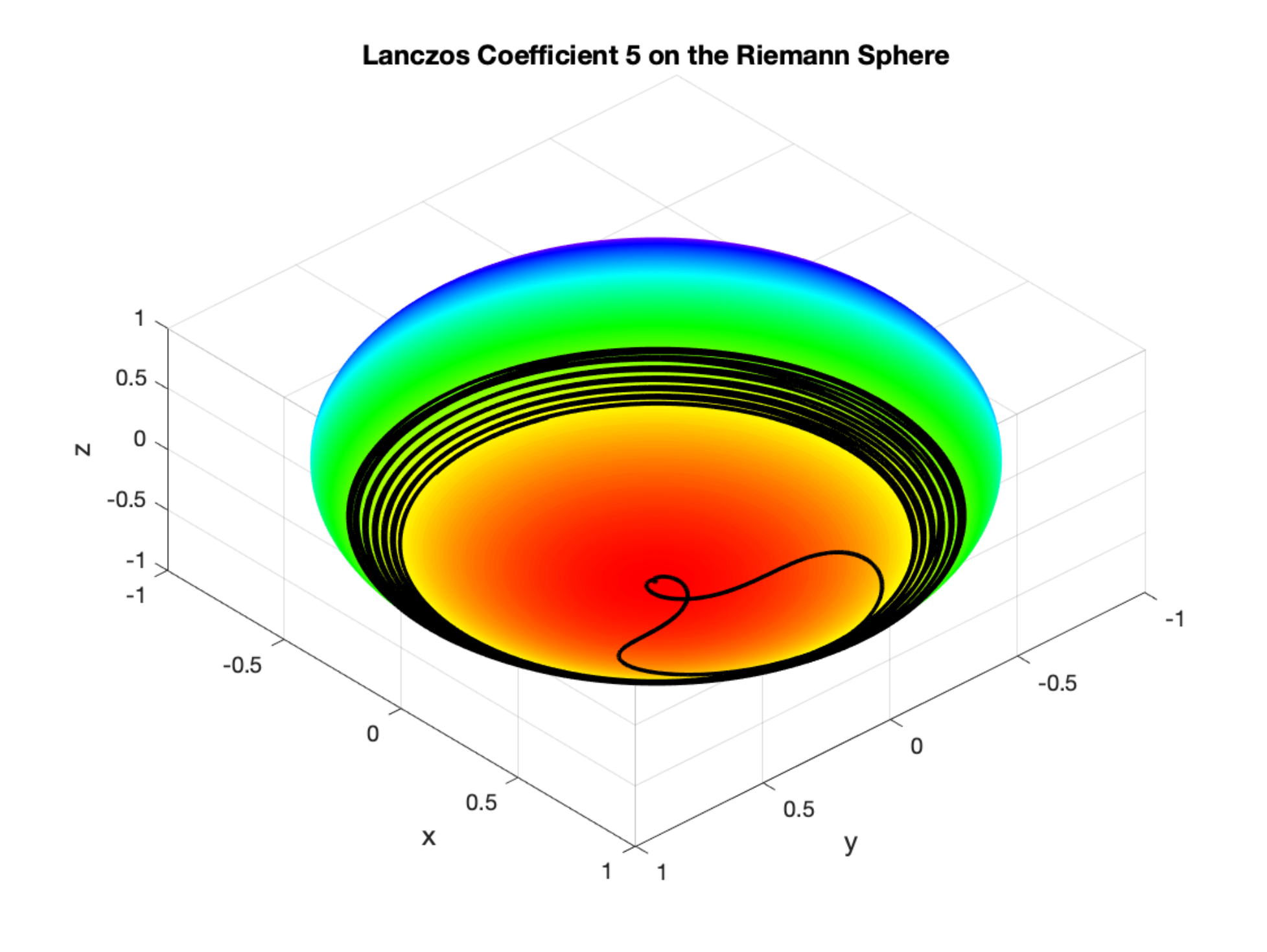}
\end{center}
\caption{The coefficient
$a_5(r)$ for $r=1+ir_y$ with $-20\pi \le r_y\le 20\pi$ plotted on the
Riemann sphere.}\label{fig:A5RS}
\end{figure}

\begin{figure}
\begin{center}
\includegraphics[width=10cm]{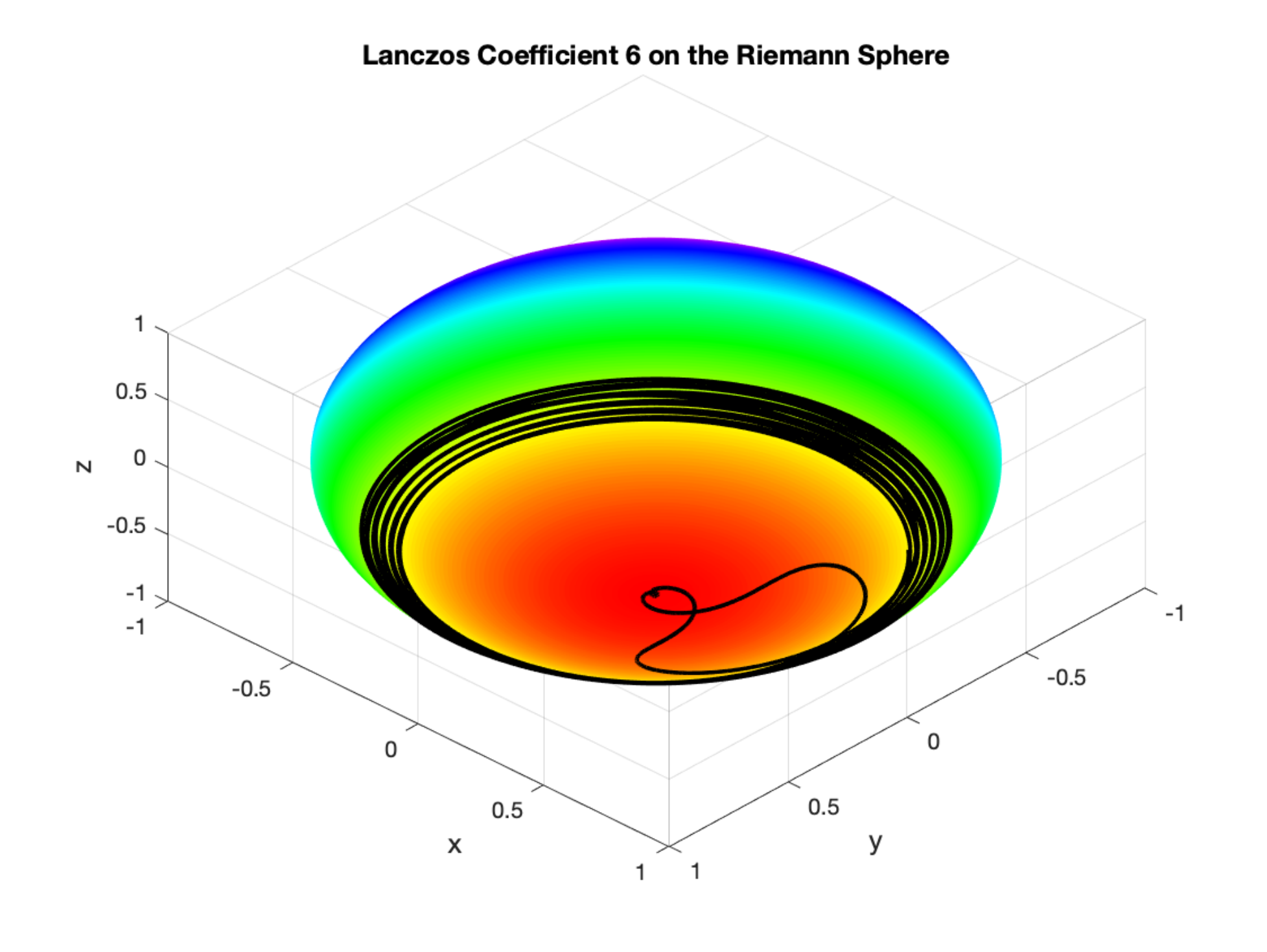}
\end{center}
\caption{The coefficient
$a_6(r)$ for $r=1+ir_y$ with $-20\pi \le r_y\le 20\pi$ plotted on the Riemann
sphere.}\label{fig:A6RS}
\end{figure}

\begin{figure}
\begin{center}
\includegraphics[width=10cm]{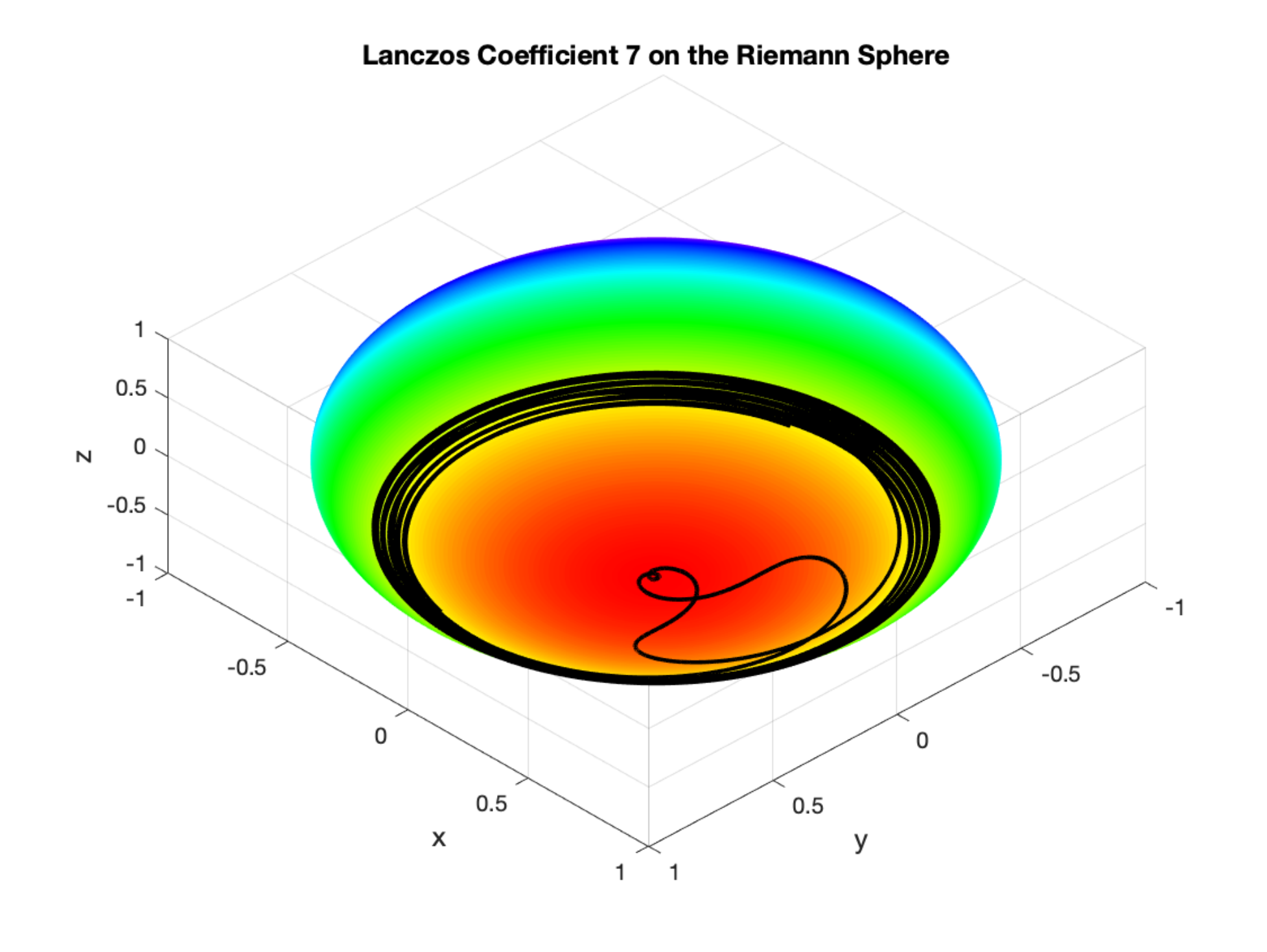}
\end{center}
\caption{The coefficient
$a_7(r)$ for $r=1+ir_y$ with $-20\pi \le r_y\le 20\pi$ plotted on the
Riemann sphere}\label{fig:A7RS}
\end{figure}

\begin{figure}
\begin{center}
\includegraphics[width=10cm]{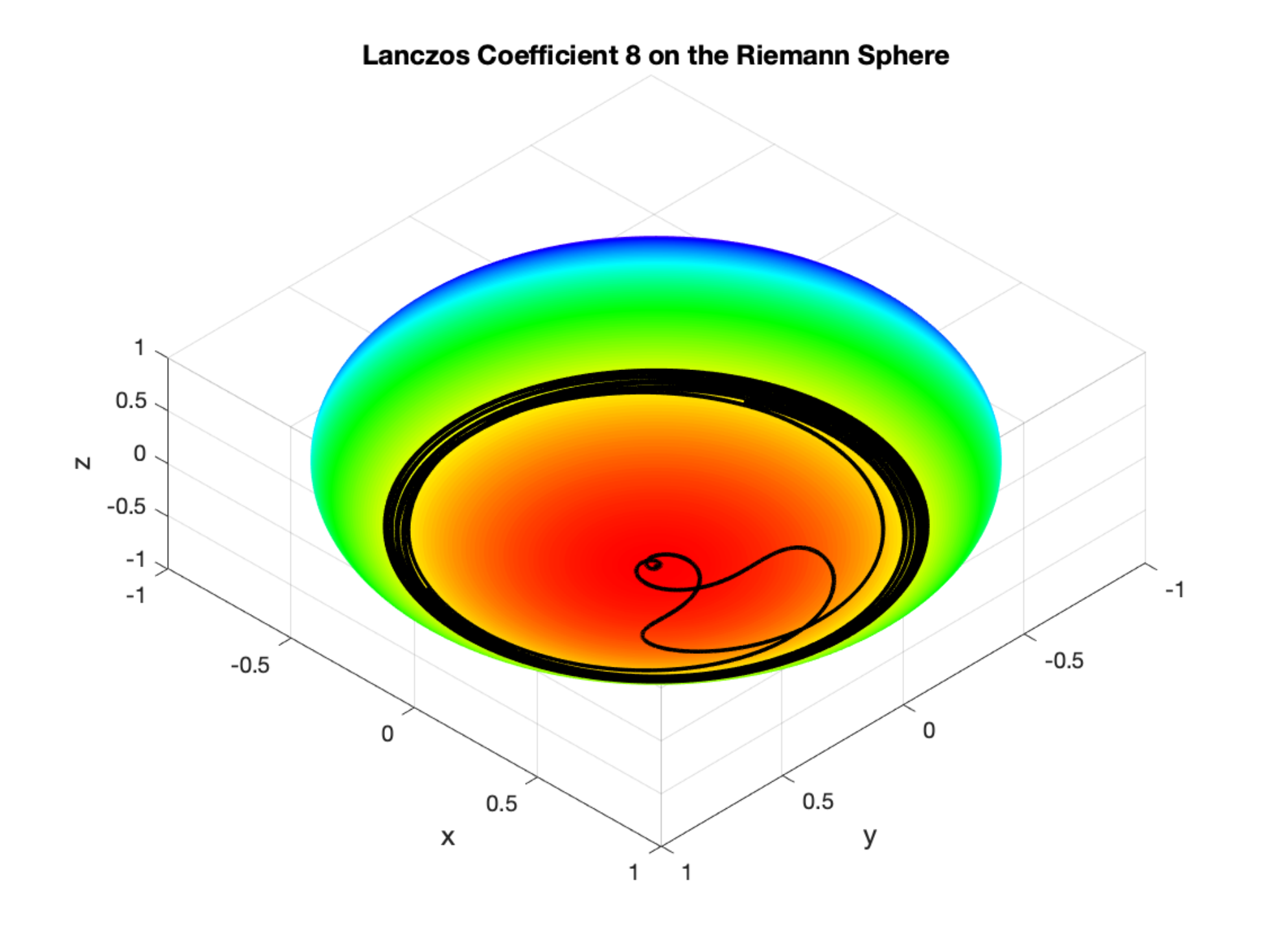}
\end{center}
\caption{The coefficient
$a_8(r)$ for $r=1+ir_y$ with $-20\pi \le r_y\le 20\pi$ plotted
on the Riemann sphere.}\label{fig:A8RS}
\end{figure}

\begin{figure}
\begin{center}
\includegraphics[width=10cm]{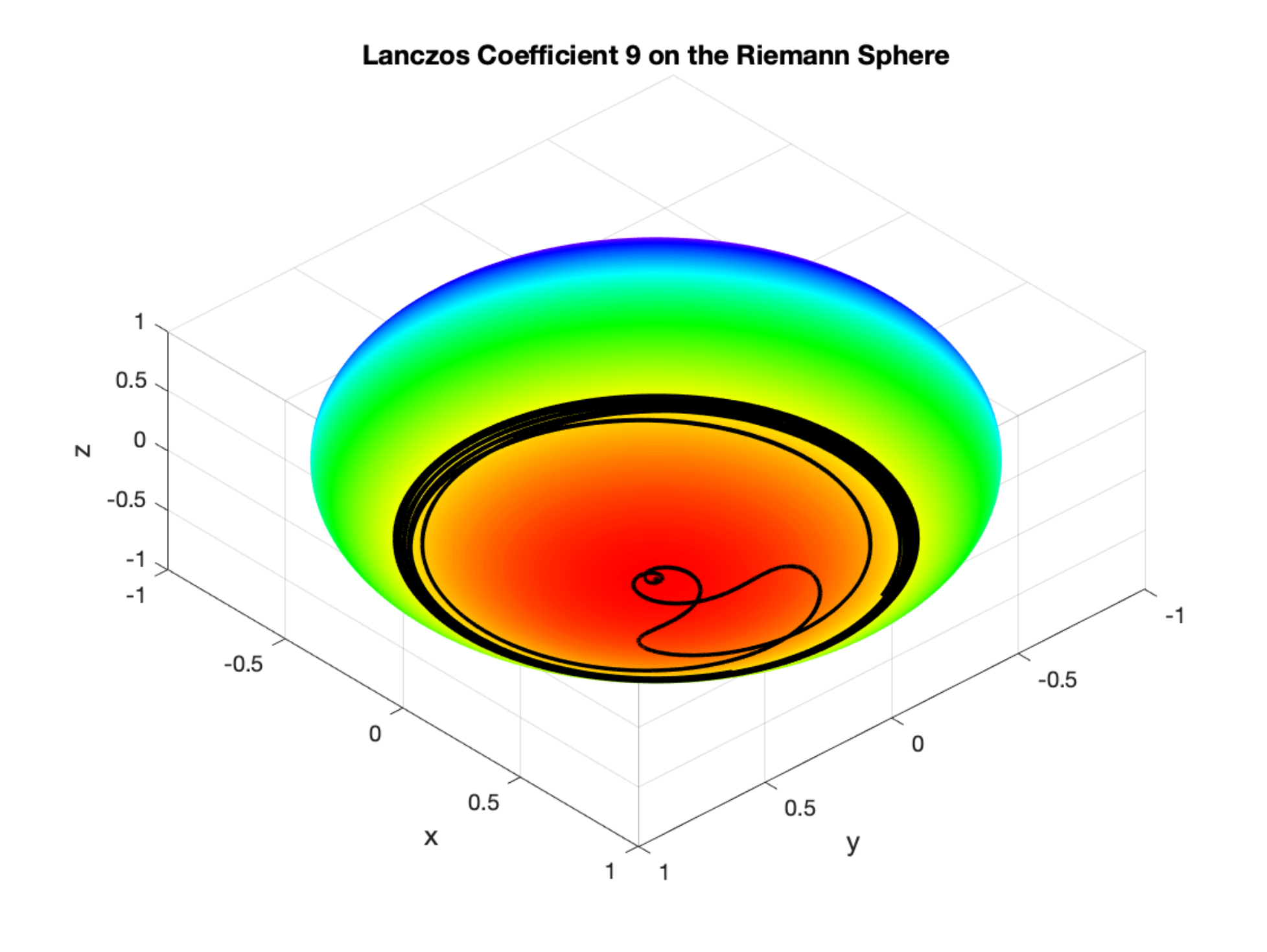}
\end{center}
\caption{The coefficient
$a_9(r)$ for $r=1+ir_y$ with $-20\pi \le r_y\le 20\pi$ plotted on the
Riemann sphere.}\label{fig:A9RS}
\end{figure}

\begin{figure}
\begin{center}
\includegraphics[width=10cm]{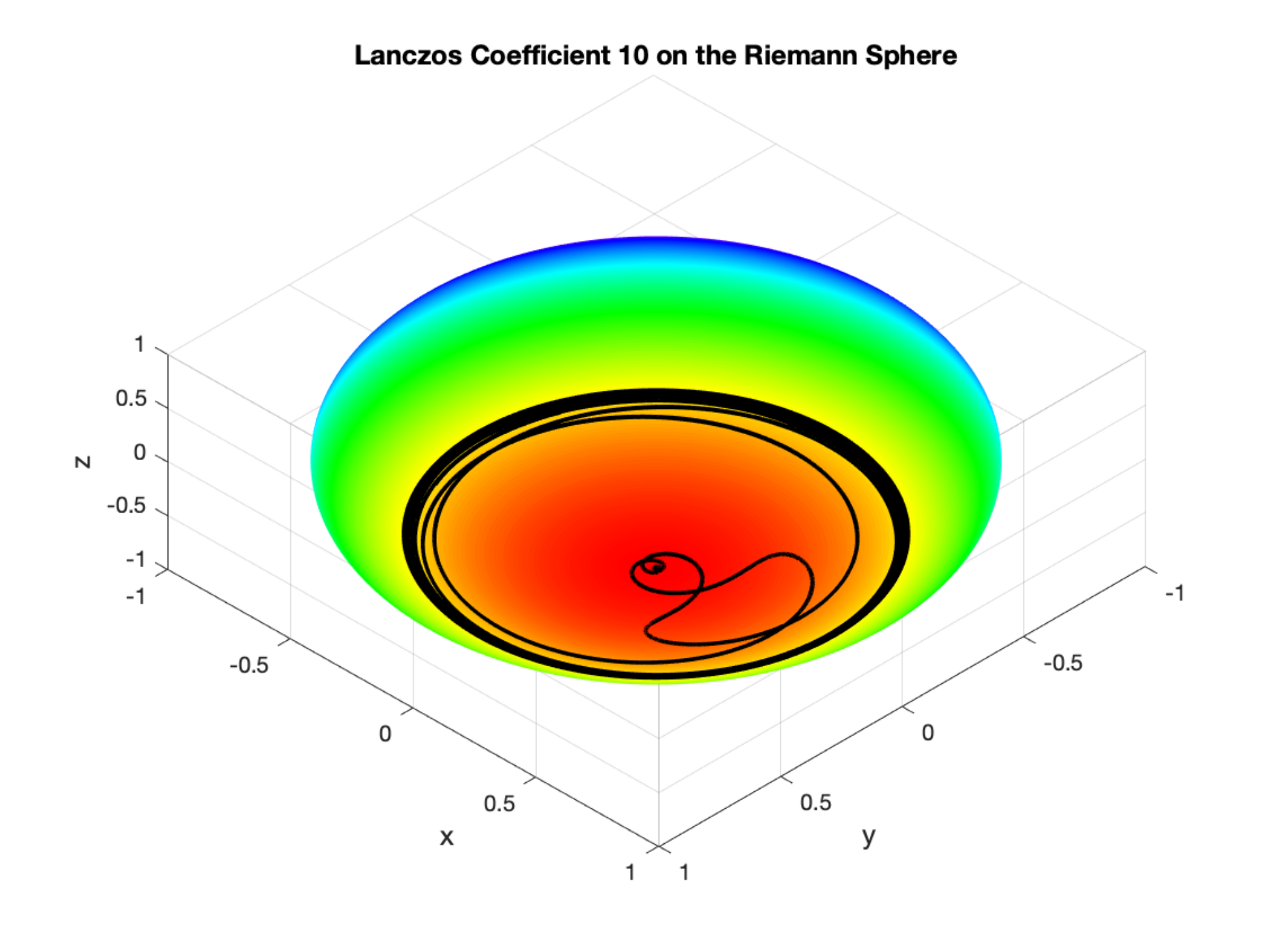}
\end{center}
\caption{The coefficient
$a_{10}(r)$ for $r=1+ir_y$ with $-20\pi \le r_y\le 20\pi$ plotted on the
Riemann sphere.}\label{fig:A10RS}
\end{figure}

\clearpage

\subsection{Magntitudes of the Coefficients}

\begin{figure}[h]
\begin{center}
\includegraphics[width=10cm]{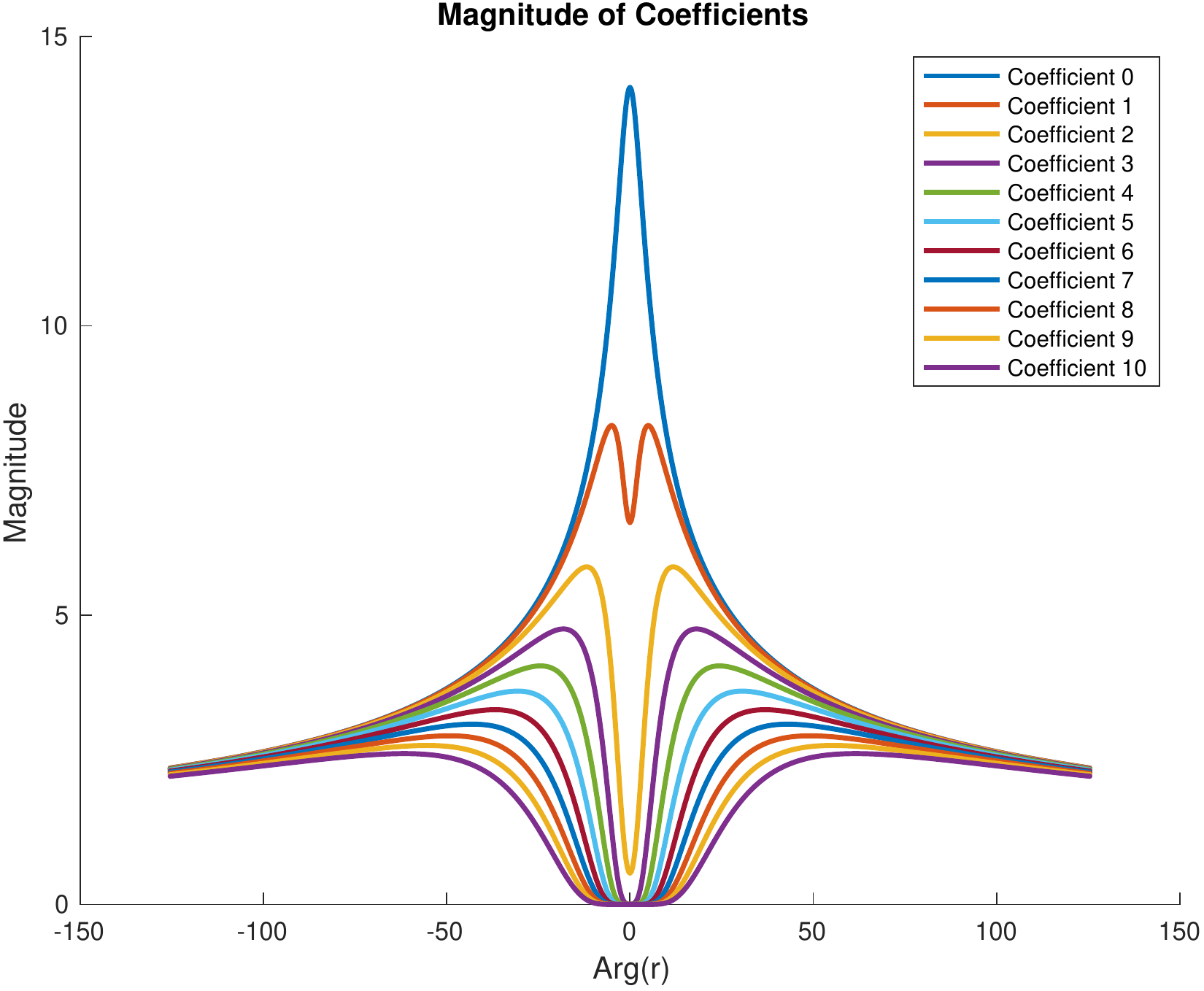}
\end{center}
\caption{A plot of the magnitudes of the first 11 coefficients for
$r=3+ir_y$ with  $-40\pi \le r_y\le 40\pi$.}\label{fig:MagnitudeR3}
\end{figure}

\begin{figure}
\begin{center}
\includegraphics[width=10cm]{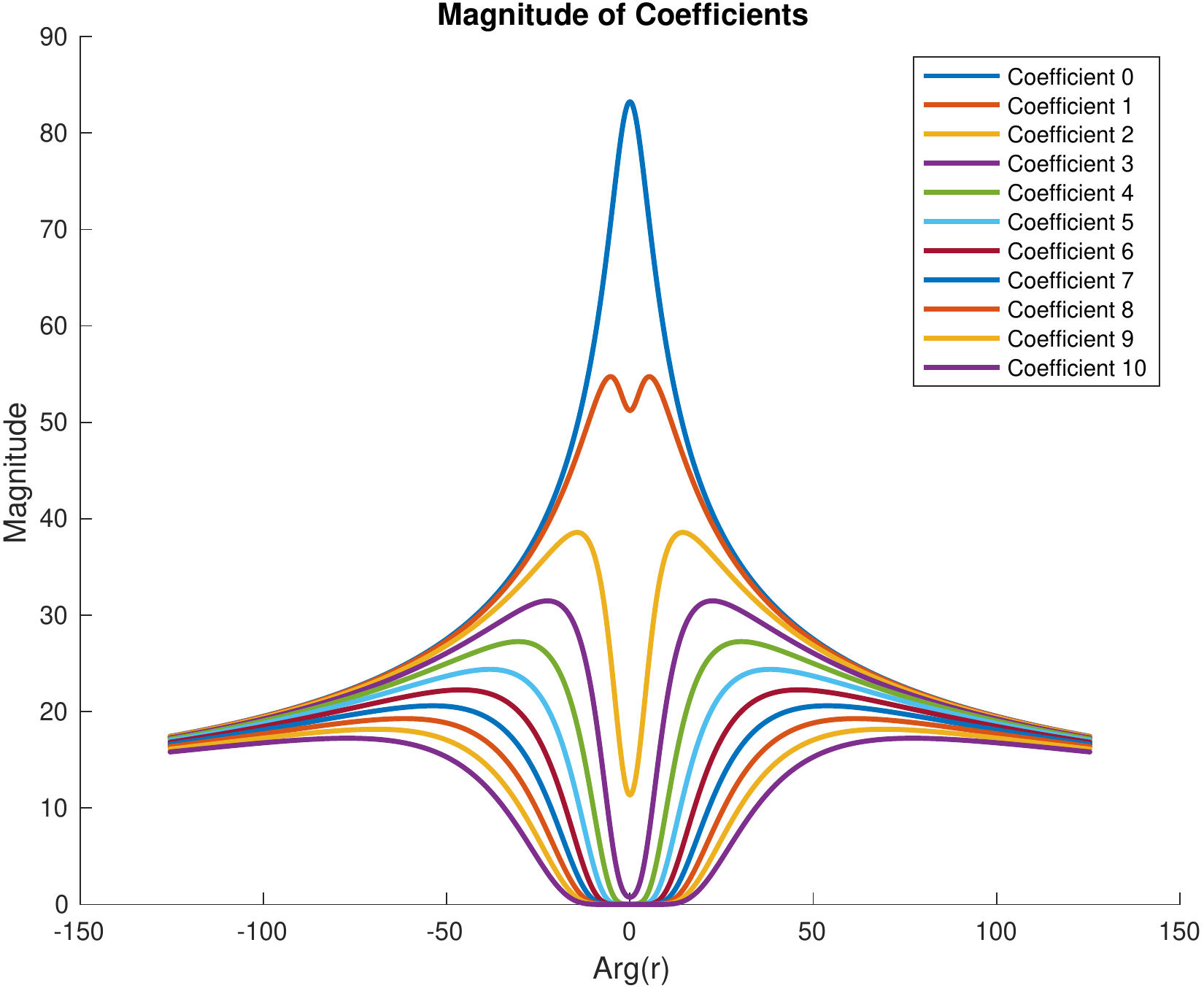}
\end{center}
\caption{A plot of the magnitudes of the first 11 coefficients for
$r=5+ir_y$ with  $-40\pi \le r_y\le 40\pi$.}\label{fig:MagnitudeR5}
\end{figure}

\begin{figure}
\begin{center}
\includegraphics[width=10cm]{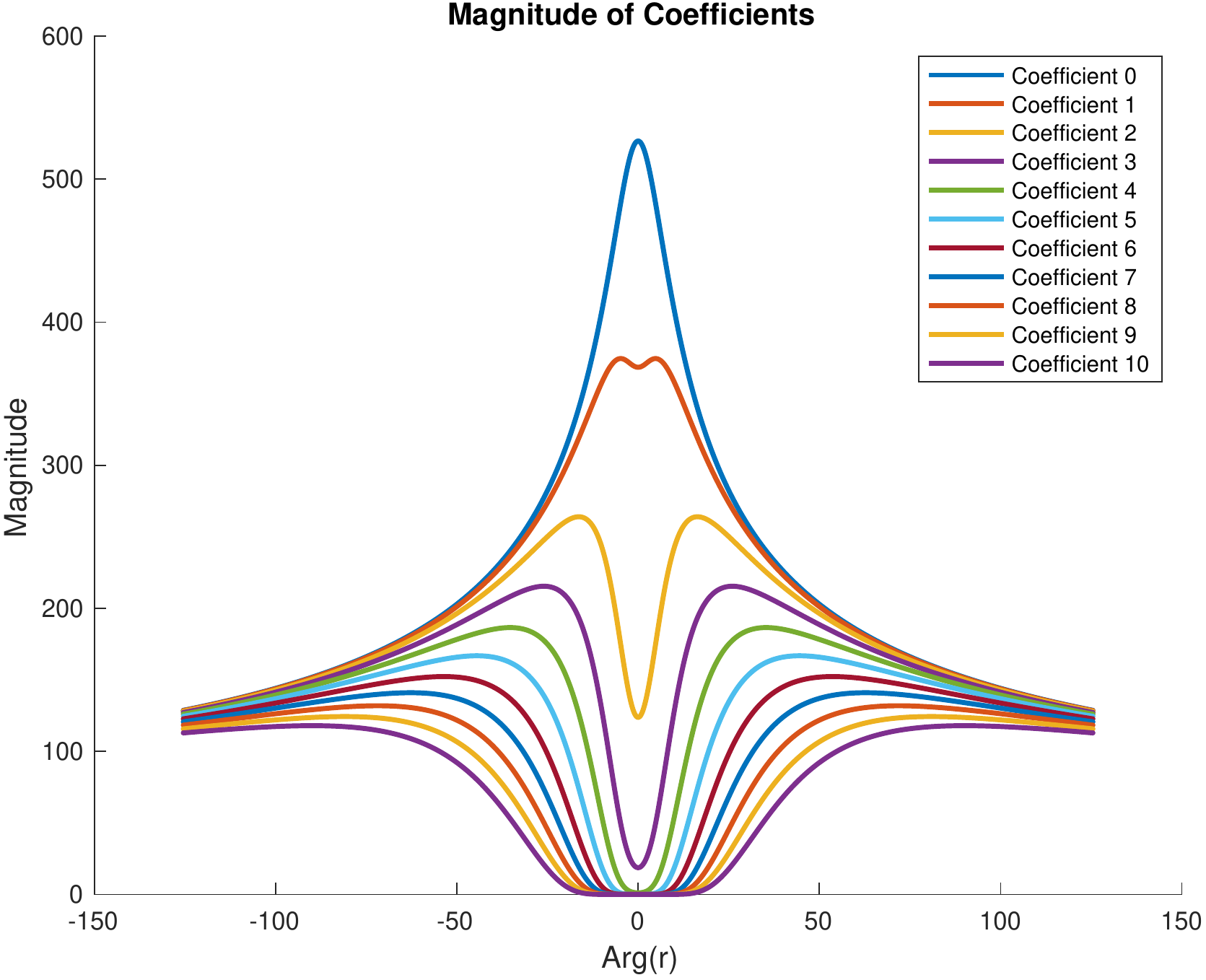}
\end{center}
\caption{A plot of the magnitudes of the first 11 coefficients for
$r=7+ir_y$ with  $-40\pi \le r_y\le 40\pi$.}\label{fig:MagnitudeR7}
\end{figure}

\begin{figure}
\begin{center}
\includegraphics[width=10cm]{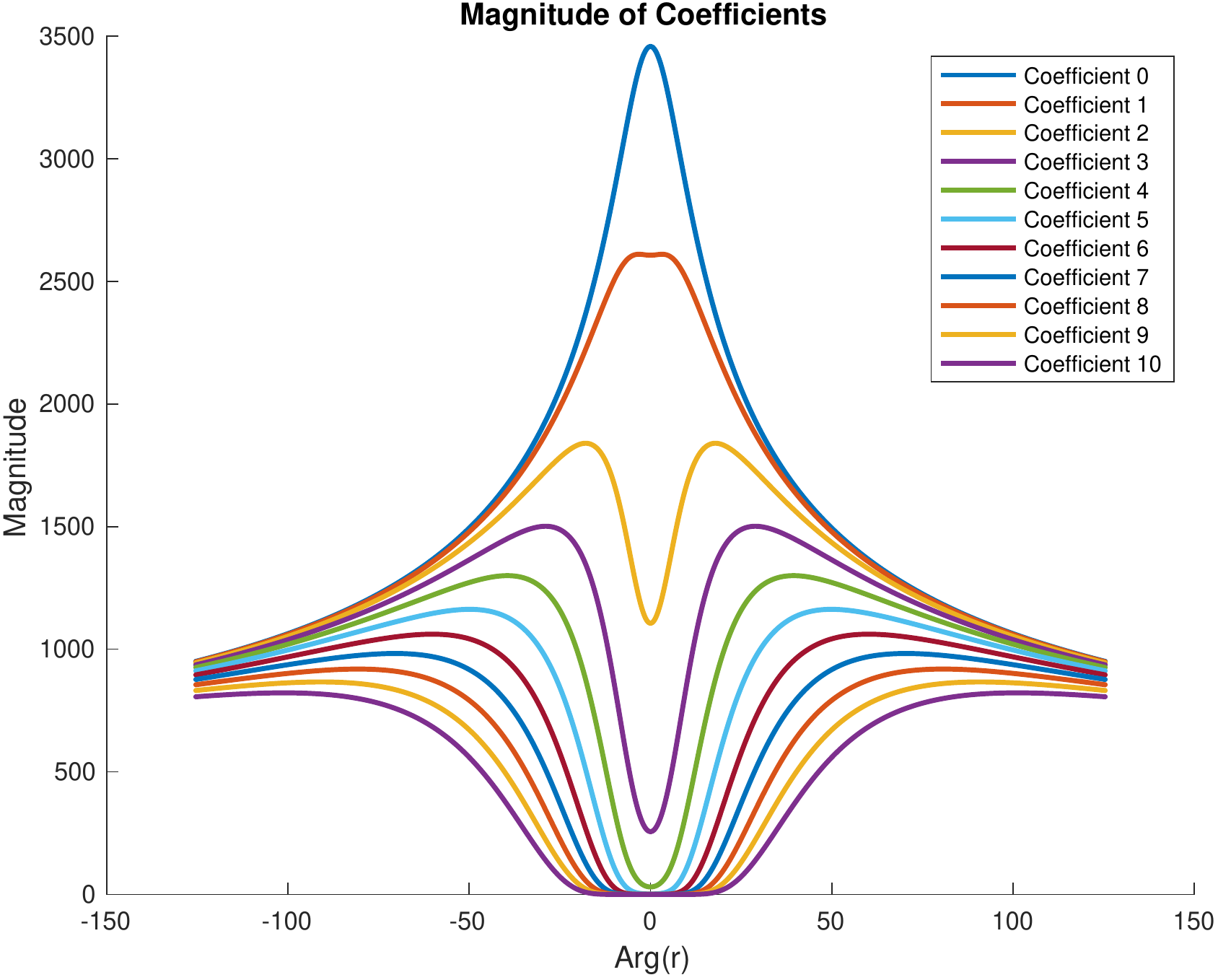}
\end{center}
\caption{A plot of the magnitudes of the first 11 coefficients for
$r=9+ir_y$ with  $-40\pi \le r_y\le 40\pi$.}\label{fig:MagnitudeR9}
\end{figure}

\begin{figure}
\begin{center}
\includegraphics[width=10cm]{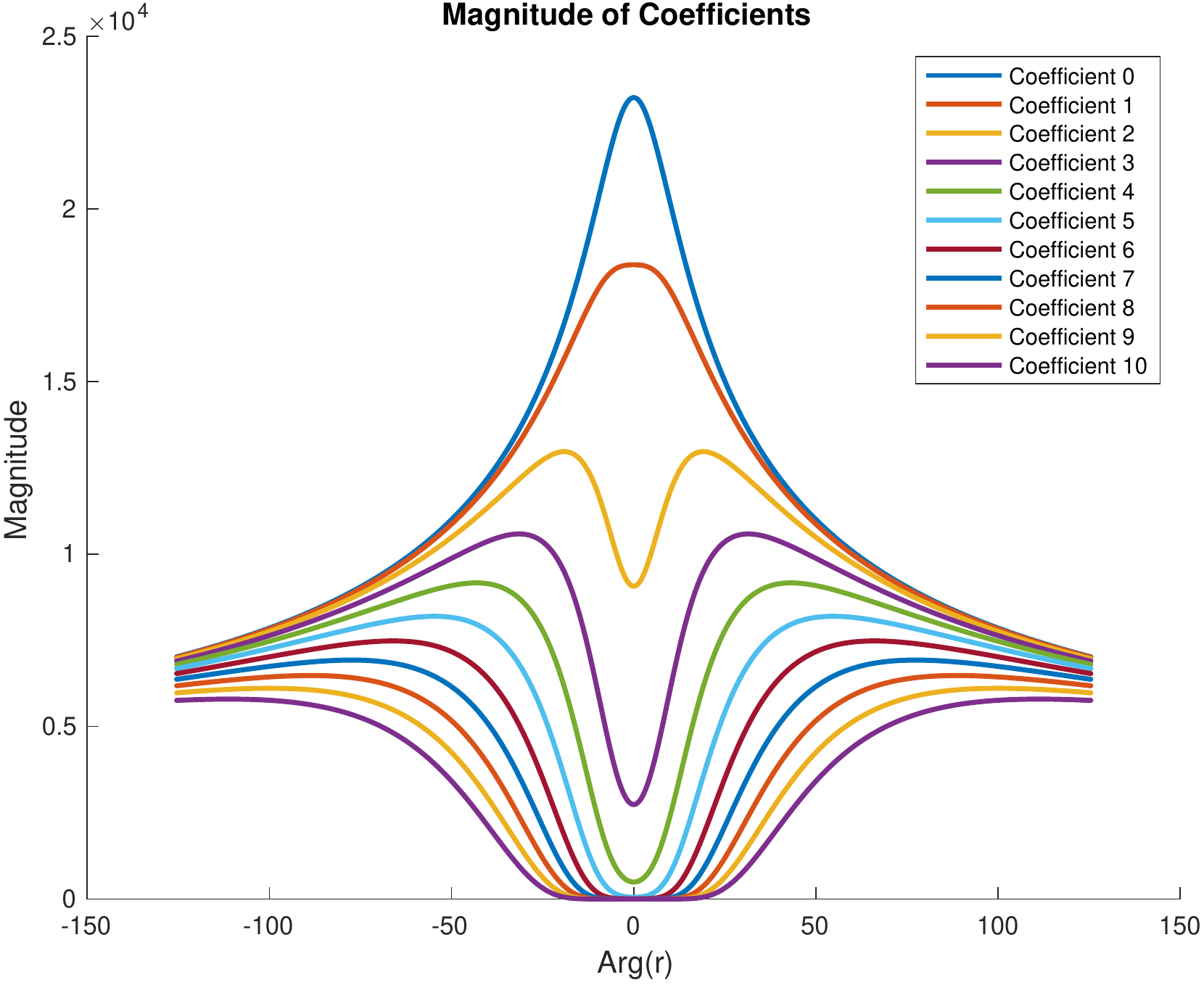}
\end{center}
\caption{A plot of the magnitudes of the first 11 coefficients for
$r=11+ir_y$ with  $-40\pi \le r_y\le 40\pi$.}\label{fig:MagnitudeR11}
\end{figure}

%% file: LanczosCoeffs.bbl
\begin{thebibliography}{}

\bibitem[\protect\citeauthoryear{Aissen}{Aissen}{1954}]{Aissen1954}
Aissen, M.~I. (1954).
\newblock {Some Remarks on Stirling's Formula}.
\newblock {\em The American Mathematical Monthly\/}~{\em 61\/}(10), 687--691.

\bibitem[\protect\citeauthoryear{Bohr and Mollerup}{Bohr and
  Mollerup}{1922}]{bohr1922}
Bohr, H. and J.~Mollerup (1922).
\newblock {\em Laereborg i Matematisk Analyse}, Volume III.
\newblock J. Gjellerup, Kopenhagen.

\bibitem[\protect\citeauthoryear{Deeba and Rodriguez}{Deeba and
  Rodriguez}{1991}]{Deeba1991}
Deeba, E.~Y. and D.~M. Rodriguez (1991).
\newblock {Stirling's Series and Bernoulli Numbers}.
\newblock {\em The American Mathematical Monthly\/}~{\em 98\/}(5), 423--326.

\bibitem[\protect\citeauthoryear{Diaconis and Freedman}{Diaconis and
  Freedman}{1986}]{diaconis1986}
Diaconis, P. and D.~Freedman (1986).
\newblock {An Elementary Proof of Stirling's Formula}.
\newblock {\em The American Mathematical Monthly\/}~{\em 93\/}(2), 123--125.

\bibitem[\protect\citeauthoryear{Dutkay, Niculescu, and Popovici}{Dutkay
  et~al.}{2013}]{Dutkay2013}
Dutkay, D.~E., C.~P. Niculescu, and F.~Popovici (2013).
\newblock {Stirling's Formula and its Extension for the Gamma Function}.
\newblock {\em The American Mathematical Monthly\/}~{\em 120\/}(8), 737--740.

\bibitem[\protect\citeauthoryear{Feller}{Feller}{1967}]{Feller1967}
Feller, W. (1967).
\newblock {A Direct Proof of Stirling's Formula}.
\newblock {\em The American Mathematical Monthly\/}~{\em 74\/}(10), 1223--1225.

\bibitem[\protect\citeauthoryear{Freitag and Busam}{Freitag and
  Busam}{2005}]{Freitag2005}
Freitag, E. and R.~Busam (2005).
\newblock {\em {Complex Analysis}}.
\newblock Springer.

\bibitem[\protect\citeauthoryear{Godfrey}{Godfrey}{2001}]{Godfrey2001}
Godfrey, P. (2001).
\newblock {A Note on the Computation of the Convergent Lanczos Complex Gamma
  Approximations}.
\newblock http://my.fit.edu/~gabdo/paulbio.html.

\bibitem[\protect\citeauthoryear{Lanczos}{Lanczos}{1964}]{lanczos1964}
Lanczos, C. (1964).
\newblock {A Precision Approximation of the Gamma Function}.
\newblock {\em Journal of the Society for Industrial and Applied Mathematics:
  Series B, Numerical Analysis\/}~{\em 1}, 86--96.

\bibitem[\protect\citeauthoryear{Lou}{Lou}{2014}]{Lou2014}
Lou, H. (2014).
\newblock {A Short Proof of Stirling's Formula}.
\newblock {\em The American Mathematical Monthly\/}~{\em 121\/}(2), 154--157.

\bibitem[\protect\citeauthoryear{Marsaglia and Marsaglia}{Marsaglia and
  Marsaglia}{1990}]{Marsaglia1990}
Marsaglia, G. and J.~C.~W. Marsaglia (1990).
\newblock {A New Derivation of Stirling's Approximation to n!}
\newblock {\em The American Mathematical Monthly\/}~{\em 97\/}(9), 826--829.

\bibitem[\protect\citeauthoryear{MATLAB}{MATLAB}{2019}]{MATLAB}
MATLAB (2019).
\newblock {\em version 9.6.0.1072779 (R2019a)}.
\newblock Natick, Massachusetts: The MathWorks Inc.

\bibitem[\protect\citeauthoryear{Mermin}{Mermin}{1984}]{Mermin1984}
Mermin, N.~D. (1984).
\newblock {Stirling's formula!}
\newblock {\em American Journal of Physics\/}~{\em 52\/}(4), 362--365.

\bibitem[\protect\citeauthoryear{Michel}{Michel}{2002}]{Michel2002}
Michel, R. (2002).
\newblock {On Stirling's Formula}.
\newblock {\em The American Mathematical Monthly\/}~{\em 109\/}(4), 388--390.

\bibitem[\protect\citeauthoryear{Namias}{Namias}{1986}]{Namias1986}
Namias, V. (1986).
\newblock {A Simple Derivation of Stirling's Asymptotic Series!}
\newblock {\em The American Mathematical Monthly\/}~{\em 93\/}(1), 25--29.

\bibitem[\protect\citeauthoryear{Neuschel}{Neuschel}{2014}]{Neuschel2014}
Neuschel, T. (2014).
\newblock {A New Proof of Stirling's Formula}.
\newblock {\em The American Mathematical Monthly\/}~{\em 121\/}(4), 350--352.

\bibitem[\protect\citeauthoryear{Patin}{Patin}{1989}]{Patin1989}
Patin, J.~M. (1989).
\newblock {A Very Short Proof of Stirling's Formula}.
\newblock {\em The American Mathematical Monthly\/}~{\em 96\/}(1), 41--42.

\bibitem[\protect\citeauthoryear{Press, Teukolsky, Vetterling, and
  Flannery}{Press et~al.}{1997}]{Press1997}
Press, W.~H., S.~A. Teukolsky, W.~T. Vetterling, and B.~P. Flannery (1997).
\newblock {\em {Numerical Recipes}}.
\newblock Cambridge University Press.

\bibitem[\protect\citeauthoryear{Pugh}{Pugh}{2004}]{Pugh2004}
Pugh, G.~R. (2004).
\newblock {\em {An Analysis of the Lanczos Gamma Approximation}}.
\newblock Ph.D. Thesis, The University of British Columbia.

\bibitem[\protect\citeauthoryear{Remmert}{Remmert}{1996}]{Remmert1996}
Remmert, R. (1996).
\newblock {Wielandt's Theorem About the $\Gamma$-Function}.
\newblock {\em The American Mathematical Monthly\/}~{\em 103\/}(3), 214--220.

\bibitem[\protect\citeauthoryear{Schmelzer and Trefethen}{Schmelzer and
  Trefethen}{2007}]{Schmelzer2007}
Schmelzer, T. and L.~N. Trefethen (2007).
\newblock {Computing the Gamma Function Using Contour Integrals and Rational
  Approximations}.
\newblock {\em SIAM Journal on Numerical Analysis\/}~{\em 45\/}(2), 558--571.

\bibitem[\protect\citeauthoryear{Schuster}{Schuster}{2001}]{Schuster2001}
Schuster, W. (2001).
\newblock {Improving Stirling's Formula}.
\newblock {\em Archiv der Mathematik\/}~{\em 77}, 170--176.

\bibitem[\protect\citeauthoryear{Spira}{Spira}{1971}]{Spira1971}
Spira, R. (1971).
\newblock {Calculation of the Gamma Function by Stirling's Formula}.
\newblock {\em Mathematics of Computation\/}~{\em 25\/}(114), 317--322.

\bibitem[\protect\citeauthoryear{Spouge}{Spouge}{1994}]{Spouge1994}
Spouge, J.~L. (1994).
\newblock {Computation of the Gamma, Digamma and Trigamma Functions}.
\newblock {\em SIAM Journal on Numerical Analysis\/}~{\em 31\/}(3), 931--944.

\end{thebibliography}
